\documentclass{amsart} 

\usepackage{amssymb} 
\usepackage{amsmath} 
\usepackage{amscd} 
\usepackage{latexsym} 
\usepackage{amsthm} 
\usepackage{verbatim} 
\usepackage{graphics} 
\usepackage{epsfig} 
\let\cal\mathcal

\def\Escr{{\cal E}} 
\def\Fscr{{\cal F}} 
\def\Gscr{{\cal G}} 
\def\Hscr{{\cal H}} 
\def\Iscr{{\cal I}} 
\def\Jscr{{\cal J}}

\def\Mscr{{\cal M}} 
\def\Nscr{{\cal N}} 
\def\Oscr{{\cal O}} 
\def\Pscr{{\cal P}}

\let\blb\mathbb 
 
\def\CC{{\blb C}}

\def \PP{{\blb P}} 
\def \ZZ{{\blb Z}} 
 
\def \NN{{\blb N}} 
\def \RR{{\blb R}} 
 
\def\LL{{\blb L}} 
\def\codim{\operatorname {codim}} 
\def\coh{\mathop{\text{\upshape{coh}}}} 
\def\coker{\operatorname {coker}} 

\def\End{\operatorname {End}} 
 
\def\Ext{\operatorname {Ext}}

\def\gldim{\operatorname {gl\,dim}}

\def\grmod{\operatorname {grmod}} 
\def\GrMod{\operatorname {GrMod}} 
\def\Hilb{\operatorname {Hilb}} 
 
\def\Hom{\operatorname {Hom}}

\def\Id{\operatorname{id}} 
\def\im{\operatorname {im}} 
\def\K0{{K_{0}(\blb P^{2}_{q})}} 
\def\ker{\operatorname {ker}}

\def\lr{\longrightarrow} 
 
\def\Mod{\operatorname{Mod}}

\def\P2q{\operatorname {{\blb P}^{2}_{q}}} 

\def\pd{\operatorname {pd}} 
\def\pdim{\operatorname {pd}}

\def\PP{\operatorname {\blb P}} 

\def\Qch{\operatorname{Qch}} 
\def\Qcoh{\operatorname {Qcoh}} 

\def\r{\rightarrow} 
 
\def\Spec{\operatorname {Spec}}

\def\Tails{\operatorname {Tails}} 
\def\tails{\operatorname {tails}} 
 
\def\Tor{\operatorname {Tor}} 

\DeclareMathOperator{\Proj}{Proj}

\DeclareMathOperator{\Aut}{Aut}

\DeclareMathOperator{\Noeth}{Noeth}

\newtheorem*{theoremA}{Theorem A} 
\newtheorem*{theoremB}{Theorem B} 
\newtheorem*{theoremC}{Theorem C} 
\newtheorem*{warning}{Warning} 
\newtheorem{lemma}{Lemma}[section] 
\newtheorem{proposition}[lemma]{Proposition} 
\newtheorem{theorem}[lemma]{Theorem} 
\newtheorem{corollary}[lemma]{Corollary} 
 
\newtheorem{convention}[lemma]{Convention}

\newtheorem{lemmas}{Lemma}[subsection] 
\newtheorem{propositions}[lemmas]{Proposition} 
 
\newtheorem{corollarys}[lemmas]{Corollary}

\theoremstyle{definition} 

\newtheorem{example}[lemma]{Example}

\newtheorem{definitions}[lemmas]{Definition}

{ 
 
\newtheorem{step}{Step}

} 

\theoremstyle{remark} 

\newtheorem{remark}[lemma]{Remark} 
\newtheorem{remarks}[lemmas]{Remark}

\newdimen\uboxsep \uboxsep=1ex 
\def\uboxn#1{\vtop to 0pt{\hrule height 0pt depth 0pt\vskip\uboxsep 
\hbox to 0pt{\hss #1\hss}\vss}} 

\def\uboxs#1{\vbox to 0pt{\vss\hbox to 0pt{\hss #1\hss} 
\vskip\uboxsep\hrule height 0pt depth 0pt}}

\numberwithin{equation}{section} 

\keywords{Weyl algebra, elliptic quantum planes, ideals, Hilbert series} 
\subjclass{Primary 16D25, 16S38, 18E30} 

\author{Koen De Naeghel and Michel Van den Bergh} 

\address{Departement WNI\\Limburgs Universitair 
Centrum\\ Universitaire Campus\\ Building D\\ 3590 
Diepenbeek\\ Belgium} 
\thanks{The second author is a director of research at the FWO} 

\email[K. De Naeghel]{koen.denaeghel@luc.ac.be} \email[M. Van den 
Bergh]{michel.vandenbergh@luc.ac.be} \date{February 29, 2004} 
\title[Ideal classes of three dimensional Artin-Schelter regular algebras]{Ideal classes of three dimensional Artin-Schelter regular algebras} 
\begin{document} 

\begin{abstract} 
  We determine the possible Hilbert functions of  graded rank one torsion 
  free modules over three dimensional Artin-Schelter regular algebras. 
  It turns out that, as in the commutative case, they are related to 
  Castelnuovo functions.  From this we obtain an intrinsic proof 
  that the space of torsion free rank one modules on a non-commutative 
  $\PP^2$ is connected. A different proof of this fact, based on 
  deformation theoretic methods and the known commutative case has 
  recently been given by Nevins and Stafford~\cite{NS}. For the Weyl 
  algebra it was proved by Wilson \cite{Wilson}. 
\end{abstract} 

\maketitle 

\tableofcontents 

\section{Introduction and main results} 
\label{ref-1-0} 
In most of this paper we work over an algebraically closed field $k$. 
Put $A=k[x,y,z]$. We view $A$ as the homogeneous coordinate ring of 
$\PP^2$. 

Let $\Hilb_n(\PP^2)$ be the Hilbert scheme of zero-dimensional 
subschemes of degree $n$ in $\PP^2$. It is well known that this is 
a smooth connected projective variety of dimension $2n$. 

Let $X\in \Hilb_n(\PP^2)$ and let $\Iscr_{X} \subset \Oscr_{\PP^{2}}$ 
be the ideal sheaf of $X$. Let $I_X$ be the graded ideal 
associated to $X$ 
\[ 
I_{X} = \Gamma_{\ast}(\PP^2,\Iscr_{X}) = 
\oplus_{l}\Gamma(\PP^{2},\Iscr_{X}(l)) 
\] 
The graded ring $A(X)=A/I_X$ is the homogeneous coordinate ring of 
$X$.  Let $h_X$ be its Hilbert function: 
\[ 
h_{X}: \NN \r \NN : m \mapsto \dim_{k} A(X)_{m} 
\] 
The function $h_X$ is of considerable interest in classical algebraic 
geometry as $h_X(m)$  gives 
the number of conditions for a plane curve of  degree $m$ to contain $X$. 
It is easy to see that $h_X(m)=n$ for $m\gg 0$, but for small values 
of $m$ the situation is more complicated (see Example 
\ref{ref-1.2-2} below).

A characterization of all possible Hilbert functions of graded ideals 
in $k[x_1,\ldots,x_n]$ was given by Macaulay in \cite{Macaulay}. 
Apparently it was Castelnuovo who first recognized the utility of the 
difference function (see \cite{Davis}) 
\[ 
s_X(m) = h_{X}(m) - h_{X}(m-1) 
\] 
Since $h_X$ is constant in high degree one has $s_X(m)=0$ for $m\gg 0$. 
It turns out that $s_X$ is a 
so-called \emph{Castelnuovo function} \cite{Davis} which 
by definition has the form 
\begin{equation} 
\label{ref-1.1-1} 
s(0)=1,s(1)=2,\ldots,s(\sigma-1)=\sigma \mbox{ and } s(\sigma-1)\ge 
s(\sigma)\ge s(\sigma+1)\ge \cdots \ge 0. 
\end{equation} 
for some integer $\sigma \geq 0$. 

It is convenient to visualize a Castelnuovo function 
using  the graph of the staircase function 
\[ 
F_{s}: \RR \r \NN: x \mapsto s({\lfloor x \rfloor}) 
\] 
and to divide the area under this graph in unit cases. We will call 
the result a \emph{Castelnuovo diagram}. The {\em weight} of a 
Castelnuovo function is the sum of its values, i.e.\ the number of 
cases in the diagram. 

  In the sequel we identify a function $f:\ZZ\r \CC$ with its 
  generating function $f(t)=\sum_n f(n) t^n$. We refer to $f(t)$ as a 
 polynomial or a series depending on whether the support of $f$ is finite 
or not.

\begin{example} 
  $s(t) = 1 + 2t + 3t^{2} + 4t^{3} + 5t^{4} + 5t^{5} + 3t^{6} + 2t^{7} 
  + t^{8} + t^{9} + t^{10}$ is a Castelnuovo polynomial of weight $28$.  The 
  corresponding diagram is 
\begin{eqnarray*} 
\includegraphics[height=2.5cm,width=5.5cm]{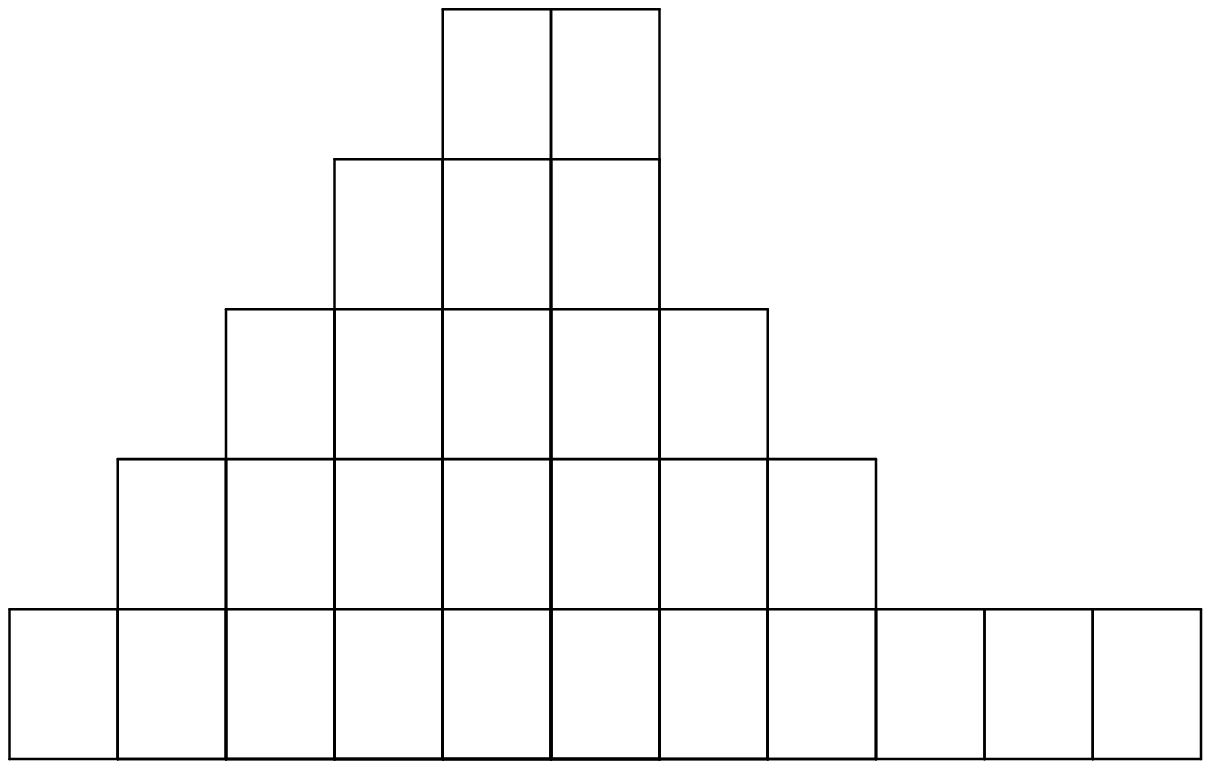} 
\end{eqnarray*} 
\end{example} 
It is known \cite{Davis, GMR, GP} that a function $h$ is of the form 
$h_X$ for $X\in \Hilb_n(\PP^2)$ if and only of $h(m)=0$ for $m<0$ and 
$h(m)-h(m-1)$ is a Castelnuovo function of weight $n$. 
\begin{example} 
\label{ref-1.2-2} Assume $n=3$. In that case there are two Castelnuovo 
  diagrams 
\[ 
\includegraphics[height=0.5cm,width=1.5cm]{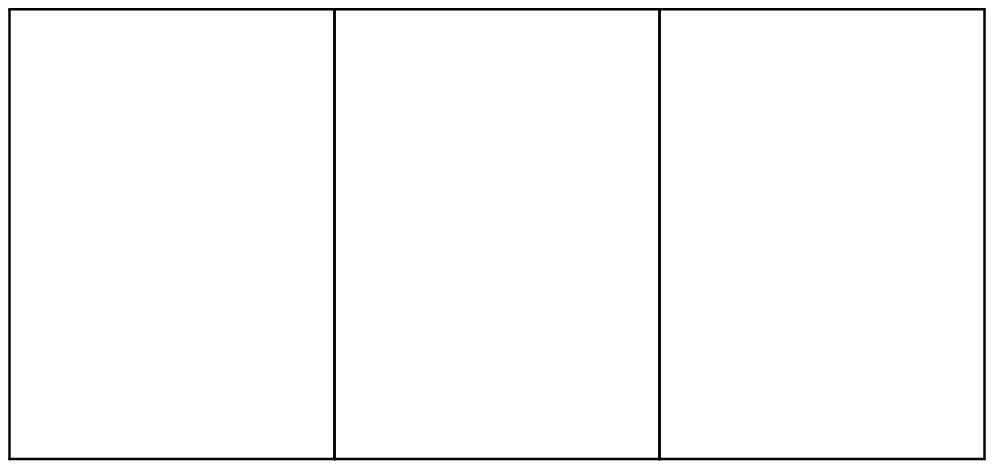}\qquad 
\includegraphics[height=1cm,width=1cm]{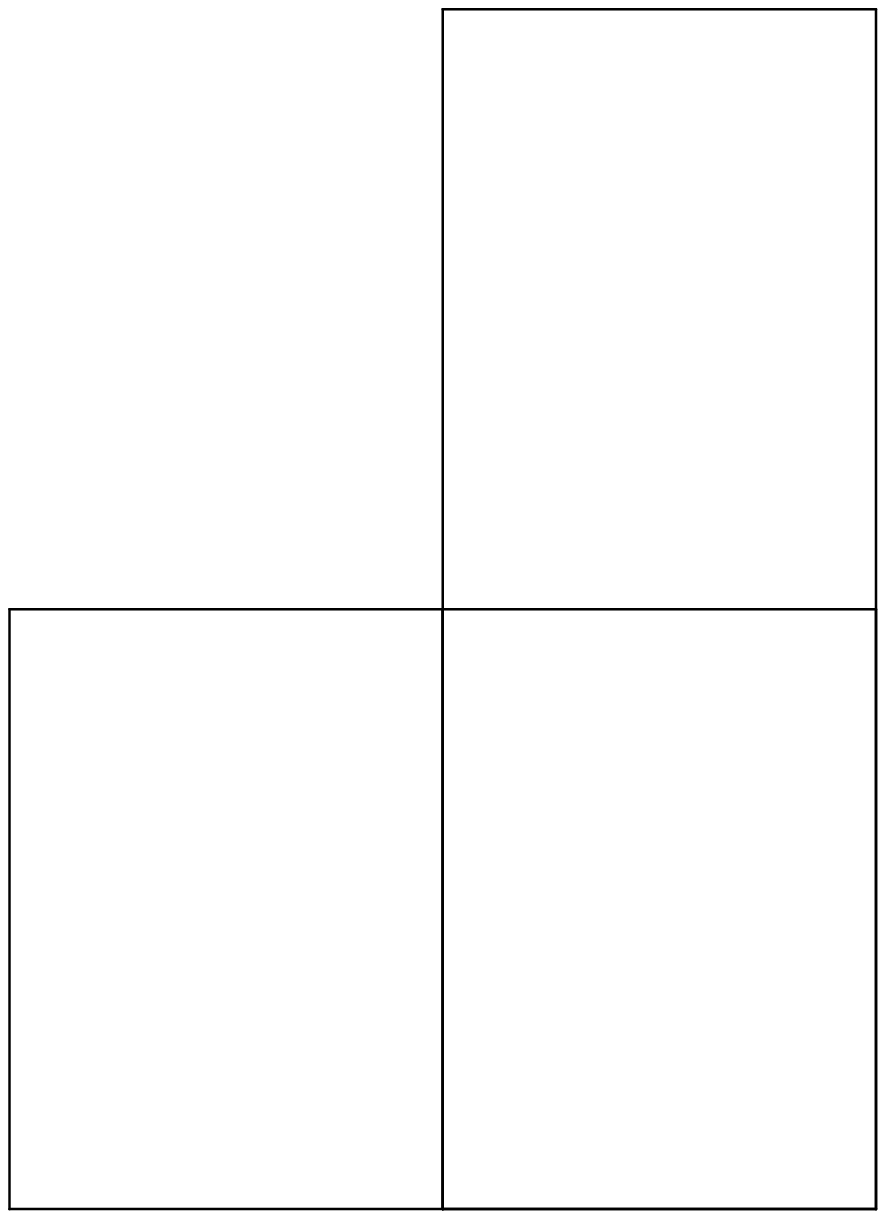}\] 
These distinguish whether the points in $X$ are collinear or not. The 
corresponding Hilbert functions are 
\[ 
1,2,3,3,3,3,\ldots \text{ and } 1,3,3,3,3,3,\ldots 
\] 
where, as expected, a difference occurs in degree one. 
\end{example} 
Our aim in this paper is to generalize the above results to the 
non-commutative deformations of $\PP^2$ which were introduced in 
\cite{AS,ATV1,ATV2,BP,OF1,OF2}. Let $A$ be a three dimensional Koszul 
Artin-Schelter regular algebra (see \S\ref{ref-3.2-16}). For the purposes of 
this introduction it suffices to say that $A$ is a non-commutative 
graded ring which is very similar to a commutative polynomial ring in 
three variables.  In particular it has the same Hilbert function and 
the same homological properties. Let $\PP^2_q$ be the corresponding 
non-commutative $\PP^2$ (see \S\ref{ref-3.1-15},\S\ref{ref-3.2-16} below). 

 The Hilbert scheme $\Hilb_n(\PP^2_q)$ was 
constructed in \cite{NS} (see also \cite{DV} for a somewhat less 
general result).  The definition of $\Hilb_n(\PP^2_q)$ is not entirely 
straightforward since in general $\PP^2_q$ will have very few 
zero-dimensional non-commutative subschemes (see \cite{smith1}), so a 
different approach is needed. It turns out that the correct generalization 
is to define $\Hilb_n(\PP^2_q)$ as the scheme parametrizing the 
torsion free graded $A$-modules $I$ of projective dimension one such that 
\[ 
h_{A}(m)-h_{I}(m)=\dim_k A_m-\dim_k I_m=n\quad\text{for}\quad m\gg 0 
\] 
(in particular $I$ has rank one as $A$-module, see \S\ref{ref-3.3-17}). 
It is 
easy to see that if $A$ is commutative then this condition singles 
out precisely the graded $A$-modules which occur as $I_X$ 
for $X\in \Hilb_n(\PP^2)$. 

The following theorem is the 
main result of this paper. 
\begin{theoremA} 
  There is a bijective correspondence between Castelnuovo polynomials 
  $s(t)$ of weight $n$ and Hilbert series $h_I(t)$ of objects in 
  $\Hilb_n(\PP^2_q)$, given by 
\begin{equation} 
\label{ref-1.2-3} 
h_I(t) = \frac{1}{(1-t)^{3}} - \frac{s(t)}{1-t} 
\end{equation} 
\end{theoremA} 
\begin{remark} By shifting the rows in a Castelnuovo diagram in such a 
  way that they are left aligned one sees that the number of diagrams 
  of a given weight is equal to the number of partitions of $n$ with 
  distinct parts. It is well-known that this is also equal to the 
  number of partitions of $n$ with odd parts \cite{Andrews}. 
\end{remark} 
\begin{remark} For the benefit of the reader we have included in 
  Appendix \ref{ref-B-68} the list of Castelnuovo diagrams of weight up to 
  six, as well as some associated data. 
\end{remark} 
From Theorem A one easily deduces that there is a 
unique maximal Hilbert series $h_{\text{max}}(t)$ and a unique minimal 
Hilbert series $h_{\text{min}}(t)$ for objects in $\Hilb_n(\PP^2_q)$. 
These correspond to the Castelnuovo diagrams 
\[ 
\parbox[c]{2.6cm}{\includegraphics[width=2.5cm]{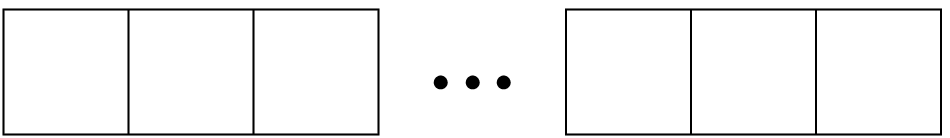}} 
\qquad\text{ and }\qquad 
\parbox[c]{2.1cm}{\includegraphics[width=2cm]{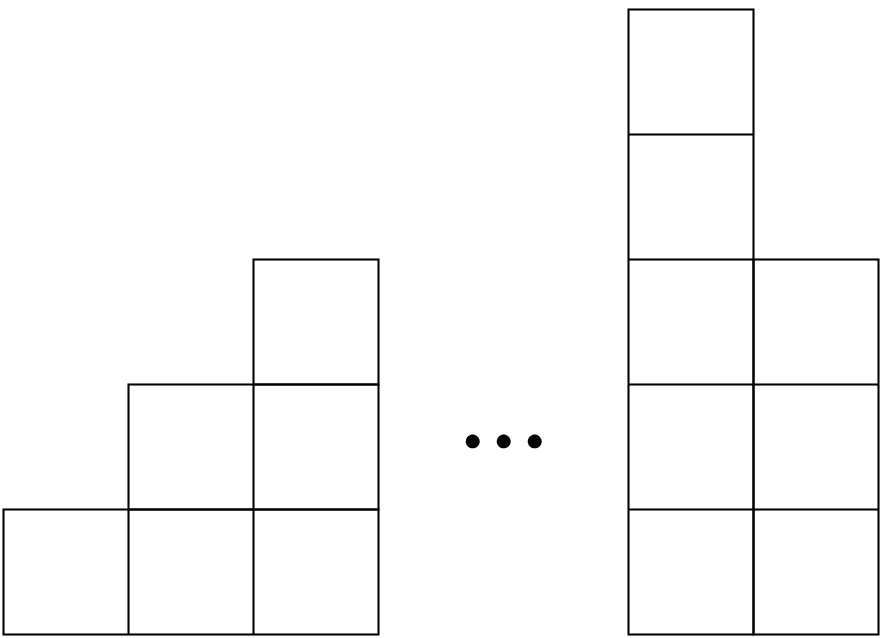}} 
\] 
We will also prove: 
\begin{theoremB}  $\Hilb_n(\PP^2_q)$ is connected. 
\end{theoremB} 
This result was recently proved for 
almost all $A$ by Nevins and Stafford \cite{NS}, using deformation 
theoretic methods and the known commutative case. In the case where 
$A$ is the homogenization of the first Weyl algebra this result was 
also proved by Wilson in \cite{Wilson}. 

We now outline our proof of Theorem B. 
For a Hilbert series $h(t)$ as in \eqref{ref-1.2-3} define 
\[ 
\Hilb_h(\PP^2_q)=\{I\in \Hilb_n(\PP^2_q)\mid h_I(t)=h(t)\} 
\] 
Clearly 
\begin{equation} 
\label{ref-1.3-4} 
\Hilb_n(\PP^2_q)=\bigcup_h \Hilb_h(\PP^2_q) 
\end{equation} 

We show below (Theorem \ref{ref-7.1-42}) that 
\eqref{ref-1.3-4} yields a stratification of $\Hilb_n(\PP^2)$ into 
non-empty smooth connected locally closed subvarieties.  In the 
commutative case this was shown by Gotzmann \cite{Gotzmann}. Our proof 
however is entirely different and seems easier. 

Furthermore 
there is a formula for $\dim \Hilb_h(\PP^2_q)$ in terms of $h$ (see 
Corollary \ref{ref-7.2.3-58} below). From that formula it follows that 
there is a \emph{unique stratum of maximal dimension} in 
\eqref{ref-1.3-4}, (which corresponds to $h=h_{\min}$). In other 
words $\Hilb_n(\PP^2_q)$ contains a dense open connected subvariety. 
This clearly implies that it is connected. 

\medskip 

To finish this introduction let us indicate how we prove Theorem 
A. Let $M$ be a torsion free graded $A$-module 
of 
projective dimension one (so we do \emph{not} require $M$ to have rank 
one). Thus 
$M$ has a minimal resolution of the form 
\begin{equation} \label{ref-1.4-5} 
0 \r \oplus_{i}A(-i)^{b_{i}} \r 
\oplus_{i}A(-i)^{a_{i}} \r M \r 0 
\end{equation} 
where $(a_{i}),(b_{i})$ are finite supported sequences of non-negative 
integers. These numbers are called the Betti numbers of $M$. They are 
related to the Hilbert series of $M$ by 
\begin{equation} 
\label{ref-1.5-6} 
h_M(t)=\frac{\sum_i (a_i-b_i)t^i}{(1-t)^3} 
\end{equation} 
So the Betti 
numbers determine the Hilbert series of $M$ but the converse is 
not true as some $a_i$ and $b_i$ may be both non-zero at the same time 
(see e.g.\ Example \ref{ref-1.7-11} below). 

Theorem A is an easy corollary of the following more refined result. 
\begin{theoremC} 
  Let $q(t)\in \ZZ[t^{-1},t]$ be a Laurent polynomial such that 
  $q_\sigma t^\sigma$ is the lowest non-zero term of $q$. Then a 
  finitely supported sequence $(a_i)$ of integers occurs among the 
  Betti numbers $(a_i),(b_i)$ of a torsion free graded $A$-module of 
  projective dimension one with Hilbert series $q(t)/(1-t)^3$ if and 
  only if 
\begin{enumerate} 
\item $a_l=0$ for $l<\sigma$. 
\item $a_\sigma=q_\sigma>0$. 
\item $\max(q_l,0)\le a_l<\sum_{i\le l} q_i$ for $l>\sigma$. 
\end{enumerate} 
\end{theoremC} 
This theorem is a natural complement to \eqref{ref-1.5-6} as it bounds the 
Betti numbers in terms of the Hilbert series. 

In Proposition \ref{ref-5.4.1-41} below we show that under suitable 
hypotheses the graded $A$-module whose existence is asserted in 
Theorem C can actually be chosen to be reflexive. This means it 
corresponds to a vector bundle on $\PP^2_q$ 
(see \S\ref{ref-3.5-23}). 

\begin{corollary} \label{ref-1.5-7} A Laurent series $h(t)=q(t)/(1-t)^3\in 
\ZZ((t))$ 
  occurs as the Hilbert series of a graded torsion free $A$-module of 
  projective dimension one if 
  and only if for some $\sigma\in \ZZ$ 
\begin{equation} 
\label{ref-1.6-8} 
\sum_{i\le l} q_i\quad 
\begin{cases} 
>0 &\text{for $l\ge \sigma$}\\ 
0&\text{for $l<\sigma$} 
\end{cases} 
\end{equation} 
I.e.\ if and only if 
\begin{equation} 
\label{ref-1.7-9} 
q(t)/(1-t)=(1-t)^2h(t)=\sum_{l\ge \sigma} p_lt^l 
\end{equation} 
with $p_l>0$ for all $l \geq \sigma$. 
\end{corollary}

In the rank one case Theorem C has the following corollary 
\begin{corollary} 
\label{ref-1.6-10} Let $h(t)=1/(1-t)^3-s(t)/(1-t)$ where $s(t)$ is a 
  Castelnuovo polynomial and let $\sigma=\max_i s_i$ (this is the same 
  $\sigma$ as in 
  \eqref{ref-1.1-1}). Then the number of minimal resolutions for an 
  object in $\Hilb_h(\PP^2_q)$ is equal to 
\[ 
\prod_{l>\sigma} [1+\min(s_{l-1}-s_l,s_{l-2}-s_{l-1})] 
\] 
This number is bigger than one if and only if there are two 
consecutive downward jumps in the coefficients of $s(t)$. 
\end{corollary} 
\begin{example} 
\label{ref-1.7-11} Assume $I\in \Hilb_n(\PP^2_q)$ has Castelnuovo diagram 
\[ 
\includegraphics[height=1cm,width=1.5cm]{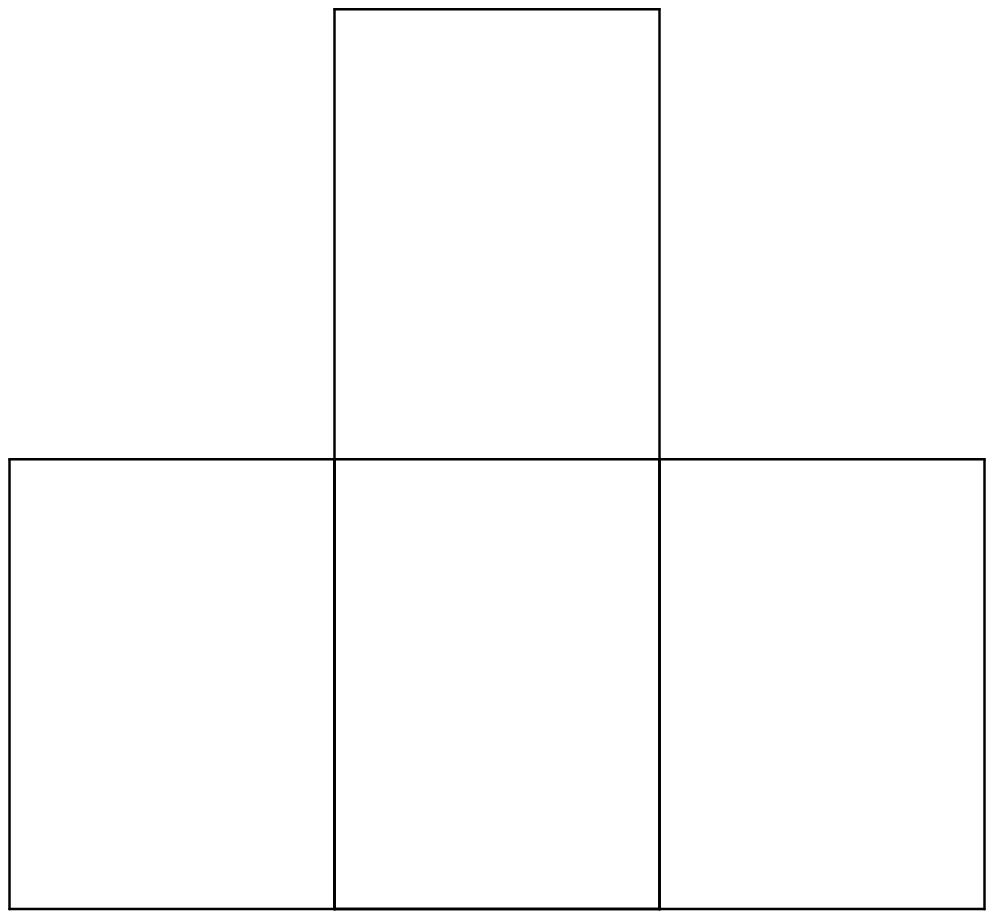}
\] 
By Corollary \ref{ref-1.6-10} we expect two different 
minimal resolutions for $I$. 
It follows from Theorem C that these are  given by 
\begin{gather} 
0  \r  A(-4)\r A(-2)^{2} \r I \r 0 \label{ref-1.8-12}\\ 
0  \r A(-3) \oplus A(-4) \r A(-2)^{2} \oplus A(-3) \r I \r 0 \label{ref-1.9-13} 
\end{gather} 
 In the commutative case 
\eqref{ref-1.8-12} corresponds to $4$ point in general position and 
\eqref{ref-1.9-13}  corresponds to a configuration of $4$ points 
among which exactly $3$ are collinear. 
\end{example}

\begin{remark} 
  Let $M$ be a torsion free graded $A$-module of projective dimension 
  one and let its Hilbert series be equal to $q(t)/(1-t)^3$. Then 
  Theorem C yields the constraint $0\le a_l< q(1)$ for $l\gg 0$ and it 
  is easy to see that $q(1)$ is equal to the rank of $M$.  Hence if 
  $M$ has rank one then there are only a finite number of 
  possibilities for its Betti numbers but this is never the case for 
  higher rank. 

     It follows that  in the case of rank $>1$ the torsion free modules $M$ of 
  projective dimension one with fixed Hilbert 
  series are not parametrized by a finite number of algebraic 
  varieties. This is to be expected as we have not imposed any 
  stability conditions on $M$. 
\end{remark}

\section{Notations and conventions} 
In this paper 
$k$ is a field which is algebraically closed except in \S\ref{ref-3-14} 
where it is arbitrary. 

Except for \S\ref{ref-7.1-48} which is about moduli spaces and 
appendix \ref{ref-A-61}, a point of a reduced scheme of finite type 
over $k$ is a closed point and we confuse such schemes with their set 
of $k$-points. 

Some results in this paper are for rank one modules and others are for 
arbitrary rank. To make the distinction clear we usually denote rank one 
modules by the letter $I$ and arbitrary rank modules by the letter 
$M$. 
\section{Preliminaries} 
\label{ref-3-14} 
In this section $k$ will be a field, not necessarily algebraically 
closed. 
\subsection{Non-commutative projective geometry} 
\label{ref-3.1-15} 
We recall some basic notions of non-commutative projective geometry. 
For more details we refer to \cite{AZ,MS,smith1,VdBSt,StICM,VdB19}. 

Let $A$ be a positively graded noetherian  $k$-algebra.  With an 
$A$-module we will mean a graded right $A$-module, and we use this 
convention for the rest of this paper.  We write $\GrMod(A)$ (resp. 
$\grmod(A)$) for the category of (resp.\ finitely generated) graded 
$A$-modules. For convenience the notations $\Hom_A(-,-)$ and 
$\Ext_A(-,-)$ will refer to $\Hom_{\GrMod(A)}(-,-)$ and 
$\Ext_{\GrMod(A)}(-,-)$. The graded $\Hom$ and $\Ext$ groups will be 
written as $\underline{\Hom}$ and $\underline{\Ext}$. 

As usual we define the non-commutative projective scheme $X=\Proj A$ of 
$A$ as the triple $(\Tails(A), \Oscr, s)$ where $\Tails A$ is the 
quotient category of $\GrMod A$ modulo the direct limits of finite 
dimensional objects, $\Oscr$ is the image of $A$ in $\Tails(A)$ and 
$s$ is the automorphism $\Mscr \mapsto \Mscr(1)$ (induced by the 
corresponding functor on $\GrMod(A)$).  We write $\Qch(X)=\Tails(A)$ 
and we let $\coh(X)$ be the noetherian objects in $\Qch(X)$. Sometimes 
we refer to the objects of $\coh(X)$ as coherent ``sheaves'' on 
$X$. 
Below it 
will be convenient to denote  objects in $\Qch(X)$ by script 
letters, like $\Mscr$. 

 We write $\pi : \GrMod(A) \r \Tails(A)$ for the 
quotient functor. The right adjoint $\omega$ of $\pi$ is given by 
$\omega\Mscr = \oplus_n\Gamma(X,\Mscr(n))$ where as usual 
$\Gamma(X,-)=\Hom(\Oscr,-)$. 
\subsection{Three dimensional Artin-Schelter regular algebras} 
\label{ref-3.2-16} 
Artin-Schelter regular algebras are non-commutative algebras which satisfy 
many of the properties of polynomial rings, therefore their associated 
projective schemes are called non-commutative projective spaces. 
\begin{definitions} \cite{AZ} 
A connected graded $k$-algebra $A$ is an {\it Artin-Schelter 
regular algebra of dimension} $d$ if it has the following properties: 
\begin{enumerate} 
\item[(i)] A has finite global dimension $d$; 
\item[(ii)] A has polynomial growth, that is, there exists positive real 
numbers $c, 
\delta$ 
such that $\dim_{k}A_{n} \leq cn^{\delta}$ for all positive integers $n$; 
\item[(iii)] A is Gorenstein, meaning there is an integer $l$ such that 
\[ 
\underline{\Ext}_{A}^{i}(k_{A},A) \cong 
\left \{ 
\begin{array}{ll} 
_{A}k(l) & \mbox{ if $i = d$,}\\ 
0 & \mbox{ otherwise.} 
\end{array} 
\right. 
\] 
where $l$ is called the {\it Gorenstein parameter} of $A$. 
\end{enumerate} 
\end{definitions} 
If $A$ is commutative then the condition (i) already implies that $A$ is 
isomorphic 
to a polynomial ring $k[x_{1}, \dots x_{n}]$ with some positive 
grading. 

There exists a complete classification for Artin-Schelter regular algebras 
of dimension three (\cite{AS} and later \cite{ATV1,ATV2,Steph1,Steph2}). 
It is known that three dimensional Artin-Schelter regular algebras 
have all expected nice homological properties. For example they are 
both left and right noetherian domains. 

In this paper we restrict ourselves to three dimensional 
Artin-Schelter regular algebras which are in addition Koszul. 
These have three generators and three defining relations in degree two. 
The minimal resolution of $k$ has the form 
\[ 
0 \r A(-3) \r A(-2)^{3} \r A(-1)^{3} \r A \r k_{A} \r 0 
\] 
hence the Hilbert series of $A$ is the same as that of the commutative 
polynomial algebra $k[x,y,z]$ with standard grading.  Such algebras 
are also referred to as {\it quantum polynomial rings in three 
  variables}.  The corresponding $\Proj A$ will be called a {\it 
  quantum projective plane} and will be denoted by $\P2q$. 

To a quantum polynomial ring in three variables $A$ we may associate 
(see \cite{ATV1}) a triple $(E, \sigma, \Oscr_E(1))$ where $E$ is either 
$\PP^{2}$ or a divisor of degree $3$ in $\PP^{2}$ (in which case we 
call $A$ {\em linear} resp. {\em elliptic}), $\Oscr_E(1) = j^{\ast} 
\Oscr_{\PP^{2}}(1)$ where $j : E \r \PP^{2}$ is the inclusion and 
$\sigma \in \Aut E$.  If $A$ is elliptic there exists, up to a scalar 
in $k$, a canonical normal element $g \in A_{3}$ and the factor ring 
$A/gA$ is isomorphic to the \emph{twisted homogeneous coordinate ring} $B = 
B(E, \Oscr_E(1), \sigma)$ (see \cite{ATV2, AVdB, AZ}). 

The fact that $A$ may be linear or elliptic presents a notational 
problem in \S\ref{ref-5.3-37}\footnote{Note that if $A$ is linear 
  then $\PP^2_q\cong \PP^2$ and we could have referred to the known 
  commutative case.} and the fact that $E$ may be non-reduced also 
presents some challenges.  We side step these problems by defining 
$C=E_{\text{red}}$ if $A$ is elliptic and letting $C$ be a $\sigma$ 
invariant line in $\PP^2$ if $A$ is linear. The geometric data 
$(E,\sigma,\Oscr_E(1))$ then restricts to geometric data 
$(C,\sigma,\Oscr_C(1))$. Denote the auto-equivalence 
$\sigma_\ast(-\otimes_C \Oscr_C(1))$ on $\Qch(C)$ by $-\otimes 
\Oscr_C(1)_\sigma$.  For $\Mscr\in \Qcoh(C)$ put $\Gamma_\ast(\Mscr)= 
\oplus_l \Gamma(C,\Mscr\otimes (\Oscr_C(1)_\sigma)^{\otimes l})$ and 
$D=B(C,\Oscr_C(1),\sigma)\overset{\text{def}}{=}\Gamma_\ast(\Oscr_C)$. 
It is easy to see that $D$ has a natural ring structure and 
$\Gamma_\ast(\Mscr)$ is a right $D$-module. Furthermore it is shown in 
\cite[Prop.\ 5.13]{ATV2} that there is a 
surjective map $A\r D=B(C,\sigma,\Mscr)$ whose kernel is generated by a 
normalizing element $h$.

 By 
analogy with the commutative case we may say that $\P2q = \Proj A$ 
contains $\Proj D$ as a ``closed'' subscheme.  Indeed it follows from 
\cite{AVdB, AZ} that the functor $\Gamma_\ast:\Qch(C)\r \GrMod(D)$ 
defines an equivalence $\Qch(C)\cong \Tails(D)$. 

The inverse of this 
equivalence and its composition with $\pi:\GrMod(D)\r \Tails(D)$ are 
both denoted by $\widetilde{(-)}$.

We define a map of non-commutative schemes  \cite{smith1} $u:C \r \P2q$ by 
\begin{align*} 
u^\ast \pi M&=(M\otimes_A D)\tilde{}\\ 
u_\ast\Mscr&=\pi(\Gamma_\ast(\Mscr)_A) 
\end{align*} 
We will call $u^\ast(\pi M)$ the \emph{restriction} of $\pi M$ to 
$C$. Clearly $u_\ast$ is an exact functor while $u^{\ast}$ is right exact. \\ 

It 
will be convenient below to let the shift functors $-(m)$ on $\coh(C)$ 
be the ones obtained from the equivalence $\coh(C)\cong \tails(D)$ and 
\emph{not} the ones coming from the embedding $C\subset \PP^2$. 
I.e. $\Mscr(n)=\Mscr\otimes (\Oscr_C(1)_\sigma)^{\otimes n}$. 
\subsection{Hilbert series} 
\label{ref-3.3-17} 
The Hilbert series of $M \in \grmod(A)$ is the Laurent power series 
\[ 
h_{M}(t) = \sum_{i = - \infty}^{+ \infty}{(\dim_{k}M_{i})t^{i}} \in\ZZ((t)). 
\] 
This definition makes sense since $A$ is right noetherian.  If 
$\Mscr\in \coh(\PP^2_q)$ then $h_\Mscr(t)=h_{\omega \Mscr}(t)$ (if $\omega 
\Mscr$ is finitely generated). 

Using a minimal resolution of $M$ we obtain the formula 
\begin{eqnarray} \label{ref-3.1-18} 
h_{M}(t)= \frac{q_{M}(t)}{ h_{A}(t)} 
\end{eqnarray} 
where $q_{M}(t)$ is an integral Laurent polynomial. 

We may write 
\begin{equation} 
\label{ref-3.2-19} 
h_{M}(t)=\frac{r}{(1-t)^3}+\frac{a}{(1-t)^2}+\frac{b}{1-t}+f(t) 
\end{equation} 
where $r,a,b\in \ZZ$ and $f(t)\in \ZZ[t^{-1},t]$. The first 
coefficient $r$ is a non-negative 
number which is called the \emph{rank} of $M$. 

If $\Mscr=\pi M$  then the numbers $r$, $a$, $b$ are determined by 
$\Mscr$. We define the rank of $\Mscr$ as the rank of $M$.

Assume $I \in \grmod(A)$ has rank one. We say that $I$ is 
\emph{normalized} if the coefficient $a$ in \eqref{ref-3.2-19} is zero. In 
that case we call $n=-b$ the \emph{invariant} of $I$.  We use the same 
terminology for objects in $\coh(\PP^2_q)$. It is easy to see that if 
$I$ has rank one then there is always an unique integer $l$ such that 
$I(l)$ is normalized. 

\begin{lemmas} Assume that $I$ is an object in 
  $\grmod(A)$. Then the following are equivalent. 
\begin{enumerate} 
\item $I$ has rank one and is normalized with invariant $n$. 
\item The Hilbert series of $I$ has the form 
\begin{equation} 
\label{ref-3.3-20} 
\frac{1}{(1-t)^3}-\frac{s(t)}{1-t} 
\end{equation} 
for a polynomial $s(t)$ with $s(1)=n$. 
\item $\dim_k A_m-\dim_k I_m=n$  for $m\gg 0$. 
\end{enumerate} 
\end{lemmas} 
\begin{proof} Easy. 
\end{proof} 
If $I$ and $s(t)$ are as in this lemma 
then we write $s_I(t)=s(t)$.  We also put $s_\Iscr(t)=s_{\omega 
  \Iscr}(t)$.

\subsection{Torsion free sheaves} 
An object in $\coh(\PP^2_q)$ or in $\grmod(A)$ is \emph{torsion} if it 
has rank zero. The corresponding notion of torsion free is defined in 
the usual way.

\begin{propositions} \label{ref-3.4.1-21} Assume that 
$\Mscr\in \coh(\PP^2_q)$ is torsion 
  free. Then $\omega \Mscr$ is finitely generated torsion free and has 
  projective dimension one. 
\end{propositions} 
\begin{proof} Assume $\Mscr=\pi M$. Without loss of generality we may 
  assume that $M$ is finitely generated and torsion free.  It follows 
  from standard localization theory that $\omega\Mscr$ is the largest 
  extension $N$ of $M$ such that $N/M$ is a union of finite length 
  modules. From this it easily follows $M\subset \omega\Mscr\subset 
  M^{\ast\ast} $ where $M^\ast=\underline{\Hom}_A(M,A)$, and hence $\omega 
\Mscr$ 
  is finitely generated and torsion free. We now replace $M$ by 
  $\omega \Mscr$. In particular $\underline{\Ext}^1_A(k,M)=0$. 

Consider a minimal resolution of $M$ 
\[ 
\cdots \r F_2\r F_1\r F_0\r M\r 0 
\] 
By applying to it the right exact functor $\underline{\Ext}_A^3(k,-)$ we 
see that 
$\underline{\Ext}^1_A(k,M)=0$ implies $F_2=0$ and hence $M$ has projective 
dimension 
one. 
\end{proof} 
\begin{corollarys} 
\label{ref-3.4.2-22} The functors $\pi$ and $\omega$ define inverse 
equivalences between the full subcategories of 
$\coh(\PP^2_q)$ and $\grmod(A)$ with objects 
\[ 
\{\text{torsion free objects in }\coh(\PP^2_q)\} 
\] 
and 
\[ 
\{\text{torsion free objects in $\grmod(A)$ of projective dimension one}\} 
\] 
\end{corollarys} 
\begin{proof} The only thing that remains to be shown is that if 
$M$ is a torsion free object in $\grmod(A)$ of projective 
dimension one then $M=\omega\pi M$. But this is clear since $\coker (M 
\r \omega\pi M)$ is finite dimensional and $\underline{\Ext}^1_A(k,M)=0$ using 
the Auslander regularity of $A$ and the fact that $\pd M=1$. 
\end{proof} 

 In case $M\in \grmod(A)$ is 
torsion free of rank one and normalized it turns out that the 
invariant of $M$ is non-negative \cite{DV,NS}. 
\subsection{Line bundles and vector bundles} 
\label{ref-3.5-23} A module $M \in \grmod(A)$ is {\em reflexive} if 
$M^{\ast \ast} = M$ where $M^{\ast} = \underline{\Hom}_{A}(M,A)$ is 
the dual of $M$.  Every reflexive module is torsion free and has 
projective dimension one.  We say that $\Mscr\in\coh(\P2q)$ is 
reflexive if this the case for $\omega\Mscr$.  If $\Mscr$ is reflexive 
 then it will be called a vector bundle. If in addition it has 
 rank one then it will be called a line bundle. 

The following criterion was proved in \cite{DV,NS} 
\begin{lemmas} 
\label{ref-3.5.1-24} 
  Assume that $A$ is an  elliptic algebra and that in the geometric data 
  $(E,\Oscr_E(1),\sigma)$ associated to $A$, $\sigma$ has infinite 
  order. If $M\in \grmod(A)$ is torsion free of projective dimension 
  one then it is reflexive if and only if its restriction $u^\ast \pi 
  M$ to the curve $C$ is a vector bundle. 
\end{lemmas} 
\section{Proof of Theorem C} 
\label{ref-5-29} 
From now on we assume that $k$ is algebraically closed. 
\subsection{Preliminaries} 
Throughout $A$ will be a quantum polynomial ring in three variables 
and $\P2q = \Proj A$ is the associated quantum projective plane. 

We will need several equivalent versions of the conditions (1-3) in the 
statement of Theorem C. One of those versions is in terms of  ``ladders''. 

For positive integers $m,n$ consider the rectangle 
\[ 
R_{m,n} = [1,m] \times [1,n] = 
\{ (\alpha, \beta) \mid 1 \leq \alpha \leq m, 1 \leq \beta \leq n \} 
\subset \ZZ^{2} 
\] 
A subset $L \subset R_{m,n}$ is called a 
{\em ladder} if 
\begin{equation*} 
\forall (\alpha,\beta) \in R_{m,n}: (\alpha,\beta) \not\in L 
\Rightarrow (\alpha+1,\beta), (\alpha,\beta-1) \not\in L 
\end{equation*} 
\begin{example} 
The ladder below is indicated with a dotted line. 

\begin{picture}(0,120) 
\put(40,0){\line(0,5){110}} 
\put(30,100){\line(5,0){100}} 
\put(36.8,0){$\vee$} 
\put(124,97.5){$>$} 
\put(27,5){$\alpha$} 
\put(137,100){$\beta$} 

\qbezier[80](45,95)(85,95)(125,95) 
\qbezier[30](45,95)(45,80)(45,65) 
\qbezier[30](45,65)(60,65)(75,65) 
\qbezier[10](75,55)(75,60)(75,65) 
\qbezier[10](75,55)(80,55)(85,55) 
\qbezier[20](85,35)(85,45)(85,55) 
\qbezier[10](85,35)(90,35)(95,35) 
\qbezier[30](95,5)(95,20)(95,35) 
\qbezier[30](95,5)(110,5)(125,5) 
\qbezier[90](125,5)(125,50)(125,95) 

\put(50,90){$\cdot$} \put(60,90){$\cdot$} \put(70,90){$\cdot$} 
\put(80,90){$\cdot$} \put(90,90){$\cdot$} \put(100,90){$\cdot$} 
\put(110,90){$\cdot$} \put(120,90){$\cdot$} 

\put(50,80){$\cdot$} \put(60,80){$\cdot$} \put(70,80){$\cdot$} 
\put(80,80){$\cdot$} \put(90,80){$\cdot$} \put(100,80){$\cdot$} 
\put(110,80){$\cdot$} \put(120,80){$\cdot$} 

\put(50,70){$\cdot$} \put(60,70){$\cdot$} \put(70,70){$\cdot$} 
\put(80,70){$\cdot$} \put(90,70){$\cdot$} \put(100,70){$\cdot$} 
\put(110,70){$\cdot$} \put(120,70){$\cdot$} 

\put(80,60){$\cdot$} \put(90,60){$\cdot$} \put(100,60){$\cdot$} 
\put(110,60){$\cdot$} \put(120,60){$\cdot$} 

\put(90,50){$\cdot$} \put(100,50){$\cdot$} \put(110,50){$\cdot$} 
\put(120,50){$\cdot$} 

\put(90,40){$\cdot$} \put(100,40){$\cdot$} \put(110,40){$\cdot$} 
\put(120,40){$\cdot$} 

\put(100,30){$\cdot$} \put(110,30){$\cdot$} \put(120,30){$\cdot$} 

\put(100,20){$\cdot$} \put(110,20){$\cdot$} \put(120,20){$\cdot$} 

\put(100,10){$\cdot$} \put(110,10){$\cdot$} \put(120,10){$\cdot$} 
\end{picture} 
\end{example} 
Let $(a_{i}),(b_{i})$ be finitely supported sequences of non-negative 
integers. We associate a sequence $S(c)$ of length $\sum_i c_i$ to 
a finitely supported sequence $(c_{i})$ as follows 
\[ 
\ldots\,,\,\underbrace{i-1,\ldots,i-1}_{c_{i-1}\text{ times }}\,,\, 
\underbrace{i,\ldots,i}_{c_i\text{ times 
}}\,,\,\underbrace{i+1,\ldots,i+1}_{c_{i+1}\text{ times }}\,,\,\ldots 
\] 
where by convention the left most non-zero entry of $S(c)$ has index 
one. 

Let $m=\sum_i a_i$, $n=\sum_i b_i$ and put $R=[1,m]\times [1,n]$. We 
associate a ladder to $(a_i),(b_i)$   as follows 
\begin{equation} 
\label{ref-5.1-30} 
L_{a,b}=\{(\alpha,\beta)\in R\mid S(a)_\alpha< S(b)_\beta\} 
\end{equation} 
\begin{lemmas} 
\label{ref-5.1.1-31} 
Let $(a_{i}),(b_{i})$ be finitely supported sequences of 
  integers and put $q_i=a_i-b_i$. The following sets of conditions are 
  equivalent. 
\begin{enumerate} 
\item Let $q_\sigma$ be the lowest non-zero $q_i$. 
\begin{enumerate} 
\item $a_l=0$ for $l<\sigma$. 
\item $a_\sigma=q_\sigma>0$. 
\item $\max(q_l,0)\le a_l<\sum_{i\le l} q_i$ for $l>\sigma$. 
\end{enumerate} 
\item Let $a_\sigma$ be the lowest non-zero $a_i$. 
\begin{enumerate} 
\item \label{ref-2a-32} The $(a_i),(b_i)$ are non-negative. 
\item \label{ref-2b-33} 
$b_{i} = 0$ for $i \leq \sigma$ 
\item \label{ref-2c-34} 
$\sum_{i\le l} b_i<  \sum_{i<l}a_{i}$ for $l > \sigma$ 
\end{enumerate} 
\item Put $m=\sum_i a_i$, $n=\sum_i b_i$. 
\begin{enumerate} 
\item The $(a_i),(b_i)$ are non-negative. 
\item $n<m$. 
\item \label{ref-3c-35} $\forall (\alpha, \beta) \in  R: \beta \geq \alpha 
- 1 \Rightarrow (\alpha, \beta)  \in L_{a,b}$. 
\end{enumerate} 
\end{enumerate} 
\end{lemmas} 
\begin{proof} The equivalence between (1) and (2) as well as the 
  equivalence 
between (2) and (3) is easy to see. We leave the details to the 
reader. 
\end{proof} 
\subsection{Proof that the conditions in Theorem C are necessary} 
We will show that the equivalent conditions given in Lemma 
\ref{ref-5.1.1-31}(2) are necessary.  The method for the proof has 
already been used in \cite{ATV2} and also by Ajitabh in 
\cite{Ajitabh}.  Assume that $M\in \grmod(A)$ is torsion free of 
projective dimension one and consider the minimal projective 
resolution of $M$. 
\begin{equation} \label{ref-5.2-36} 
0 \r \oplus_{i}A(-i)^{b_{i}} \r \oplus_{i}A(-i)^{a_{i}} \r M \r 0 
\end{equation} 
There is nothing to prove for (\ref{ref-2a-32}) so we discuss 
(\ref{ref-2b-33})(\ref{ref-2c-34}). Since \eqref{ref-5.2-36} is 
a minimal resolution, it contains for all integers $l$ a subcomplex of 
the form 
\[ 
\oplus_{i \le l}  A(-i)^{b_i} \xrightarrow{\phi_l} \oplus_{i< l} A(-i)^{a_i} 
\] 
The fact that $\phi_l$ must be injective implies 
\[ 
\sum_{i \le l}{b_i} \leq \sum_{i < l}{a_i} 
\] 
In particular, if we take $l = \sigma$ this already shows that 
$b_{i} = 0$ for $i \leq \sigma$ which proves (\ref{ref-2b-33}). 
Finally, to prove (\ref{ref-2c-34}), assume that there is some $l > \sigma$ 
such that $\sum_{i \le l}{b_i} = \sum_{i < l}{a_i}$. 
This means that $\coker \phi_l$ is torsion and different from zero. 
Note that $\oplus_{i < l} A(-i)^{a_i}$ is not zero since $l > \sigma$. 
We have a map 
\[ 
\coker \phi_l\r M 
\] 
which must be zero since $M$ is assumed to be torsion free. 
But this implies that $ \oplus_{i < l} A(-i)^{a_i}\r M$ is the zero map, 
which is obviously impossible given the minimality of our chosen resolution 
\eqref{ref-5.2-36}. 
Thus we obtain that 
\[ 
\sum_{i \le l}{b_i} < \sum_{i < l}{a_i} 
\] 
which completes the proof. 
\subsection{Proof that the conditions in Theorem C are sufficient.} 
\label{ref-5.3-37} 
In this section the notations and conventions are as in 
\S\ref{ref-3.1-15},\S\ref{ref-3.2-16}. 

We will assume the equivalent conditions given in Lemma 
\ref{ref-5.1.1-31}(3) hold. Thus we fix finitely supported sequences 
$(a_i)$, $(b_i)$ of non-negative integers such that $n=\sum_i 
b_i<m=\sum_i a_i$ and we assume in addition that the ladder condition 
(\ref{ref-3c-35}) is true.

Our proof of the converse of Theorem C is a suitably  adapted version of 
\cite[p468]{CGO}. It is based on a series of observations, the first one 
of which is the next lemma. 
\begin{lemmas} If $M\in \grmod(A)$ has a resolution (not necessarily minimal) 
\[ 
0\r \oplus_i A(-i)^{b_i}\xrightarrow{\phi} \oplus_i A(-i)^{a_i}\r M\r 0 
\] 
such that the restriction 
\[ 
u^\ast(\pi \phi): \oplus_i \Oscr_{C}(-i)^{b_i}\r \oplus_i 
\Oscr_{C}(-i)^{a_i} 
\] 
has maximal rank at every point in $C$ then $M$ is torsion free. 
\end{lemmas} 
\begin{proof} Assume that $M$ is not torsion free and that $u^\ast(\pi \phi)$ 
  has the stated property. This means that $u^\ast(\pi \phi)$ is an 
  injective map whose cokernel $u^\ast \pi M$ is a vector bundle on~$C$. 

Let $T$ be the torsion submodule of $M$. 
  Note first that $M$ cannot have a submodule of GK-dimension $\le 1$ 
  as $\Ext^1_A(-,A)$ is zero on modules of GK-dimension $\le 1$ 
  \cite{ATV2}. Hence $T$ has pure GK-dimension two. 

  If $T$ contains $h$-torsion then $\Tor_1^A(D,M)$ is not zero and in 
  fact has GK-dimension two.  Thus $u^\ast(\pi \phi)$ is not 
  injective, yielding a contradiction. 

Assume now that $T$ is $h$-torsion free. In that case $T/Th$ is a 
submodule of GK-dimension one of $M/Mh$. And hence $u^\ast \pi T$ is a 
submodule of dimension zero of $u^\ast \pi M$ which is again a 
contradiction. 
\end{proof} 
Now note that the map 
\begin{equation} 
\label{ref-5.3-38} 
\Hom_A(\oplus_i A(-i)^{b_i}, \oplus_i A(-i)^{a_i})\r 
\Hom_C(\oplus_i \Oscr_{C}(-i)^{b_i},\oplus_i 
\Oscr_{C}(-i)^{a_i}):\phi\mapsto u^\ast(\pi\phi) 
\end{equation} 
is surjective.  Let $H$ be the linear subspace of $\Hom_C(\oplus_i 
\Oscr_{C}(-i)^{b_i},\oplus_i \Oscr_{C}(-i)^{a_i})$ whose elements are 
such that the projections on $\Hom_C( \Oscr_{C}(-i)^{b_i}, 
\Oscr_{C}(-i)^{a_i})$ are zero for all $i$.  If we can find $N\in H$ 
of maximal rank in every point then an arbitrary lifting of $N$ under 
\eqref{ref-5.3-38} yields a torsion free $A$-module with Betti numbers 
$(a_i)$, $(b_i)$. 

The elements of $H$ are given by matrices 
$(h_{\alpha\beta})_{\alpha\beta}$ for $(\alpha,\beta)\in L_{a,b}$ 
where $L_{a,b}$ is as in \eqref{ref-5.1-30} and 
where the 
$h_{\alpha\beta}$ are elements of  suitable non-zero 
$\Hom_C(\Oscr_C(-i),\Oscr_C(-j))$. We will look for $N$ in the linear 
subspace ${}^0\!H$ of $H$ given by those matrices where 
$h_{\alpha\beta}=0$ for $\beta\neq \alpha,\alpha-1$. 

\medskip 

To construct find $N$ we use the next 
observation. 
\begin{lemmas} 
\label{ref-5.3.2-39} For $p\in C$  and $N\in {}^0\!H$ let $N_p$ be the 
restriction of $N$ to 
  $p$ and write 
\[ 
{}^0\!H_p=\{N\in {}^0\!H\mid \text{$N_p$ has  non-maximal rank}\} 
\] 
If 
\begin{equation} 
\label{ref-5.4-40} 
\codim_{{}^0\!H} {}^0\!H_p\ge 2\text{ for all }p\in C 
\end{equation} 
 then there exists an $N$ 
in ${}^0\!H$ which has maximal rank everywhere. 
\end{lemmas} 
\begin{proof} 
Assume that \eqref{ref-5.4-40} holds. 
  Since $({}^0\!H_p)_p$ is a one-dimensional family of subvarieties of 
  codimension $\ge 2$ in ${}^0\!H$ it is intuitively clear that their union 
  cannot be the whole of ${}^0\!H$, proving the lemma. 

  To make this idea precise let $\Escr_1$, $\Escr_0$ be the pullbacks 
  of the vector bundles $\oplus_i \Oscr_{C}(-i)^{b_i}$, $\oplus_i 
  \Oscr_{C}(-i)^{a_i}$ to ${}^0\!H\times C$ and let $\Nscr:\Escr_1\r 
  \Escr_0$ be the vector bundle map which is equal to $N_p$ in the 
  point $(N,p)\in {}^0\!H\times C$. Let ${}^0\!\Hscr\subset H\times 
  C$ be the locus of points $x$ in ${}^0\!H\times C$ where $\Nscr_x$ 
  has non-maximal rank. It is well-known and easy to see that 
  ${}^0\!\Hscr$ is closed in ${}^0\!H\times C$. A more down to earth 
  description of ${}^0\!\Hscr$ is 
\[ 
 {}^0\!\Hscr=\{(N,p)\in {}^0\!H\times C\mid \text{$N_p$ has 
 non-maximal rank}\} 
\] 
By considering the fibers of the projection ${}^0\!H\times C\r C$ we see that 
${}^0\!\Hscr$ has codimension $\le 2$ in ${}^0\!H\times C$. Hence its 
projection on ${}^0\!H$, which is $\bigcup_p {}^0\!H_p$, has codimension 
$\ge 1$. 
\end{proof} 
Fix a point $p\in C$ and fix basis elements for the one-dimensional 
vector spaces $\Oscr_C(-i)_p$. Let $\LL$ be the vector spaces 
associated to the ladder $L_{a,b}$ (see \eqref{ref-5.1-30}) as follows 
\[ 
\LL = \{ A \in M_{m \times n}(k) \mid A_{\alpha \beta} = 0 \mbox{ for } 
(\alpha,\beta) \not\in L \} 
\] 
and 
let ${}^0\LL$ be the subspace defined by $A_{\alpha\beta}=0$ for 
$\beta\neq \alpha,\alpha-1$. 
Then there is a surjective linear map 
\[ 
\phi_p:{}^0\!H\mapsto {}^0\LL:N\mapsto N_p 
\] 
Let $V$ be the matrices of non-maximal rank in ${}^0\LL$. We have 
\[ 
{}^0\!H_p=\phi_p^{-1}(V) 
\] 
Now by looking at the two topmost $n\times n$-submatrices 
we see 
that for a matrix in ${}^0\LL$ to not have  maximal rank both the 
diagonals $\beta=\alpha$ and $\beta=\alpha-1$ must contain a zero 
(this is not sufficient). Using condition \ref{ref-5.1.1-31}(3c) we 
see that $V$ has codimension $\ge 2$ and so the 
same holds for ${}^0\!H_p$. This means we are done. 
\begin{remarks} It is easy to see that the actual torsion free module 
  constructed in this section is the direct sum of a free module and a 
  module of rank one. 
\end{remarks} 
\subsection{A refinement} 

\begin{propositions} 
\label{ref-5.4.1-41} Assume that $A$ is a elliptic and that in 
the geometric data 
$(E,\Oscr_E(1),\sigma)$ associated to $A$, $\sigma$ has infinite 
order. Then the graded $A$-module whose existence is asserted in 
Theorem C can  be chosen to be reflexive. 
\end{propositions} 
\begin{proof} 
  The modules that are constructed in \S\ref{ref-5.3-37} satisfy the 
  criterion given in Lemma \ref{ref-3.5.1-24}, hence they are 
  reflexive. 
\end{proof}

\section{Proof of other properties of Hilbert series} 
\begin{proof}[Proof of Corollary \ref{ref-1.5-7}] 
It is easy to see that 
the conditions (1-3) in Theorem C have a solution for $(a_i)$ if 
and only if \eqref{ref-1.6-8} is true. The equivalence of \eqref{ref-1.6-8} and 
\eqref{ref-1.7-9} is clear. 
\end{proof} 
\begin{proof}[Proof of Theorem A] 
Let $h(t)$ is a Hilbert series of the form \eqref{ref-3.3-20}. Thus 
$h(t)=q(t)/(1-t)^3$ where $q(t)=1-(1-t)^2s(t)$ 
and hence 
$q(t)/(1-t)=1/(1-t)-(1-t)s(t)$. Thus \eqref{ref-1.7-9} is equivalent to 
$(1-t)s(t)$ being of the form 
\[ 
(1-t)s(t)=1+t+t^2+\cdots+t^{\sigma-1}+d_\sigma t^\sigma + \cdots 
\] 
where $d_i\le 0$ for $i\ge \sigma$. Multiplying by 
$1/(1-t)=1+t+t^2+\cdots$ shows that this is equivalent to $s(t)$ being 
a Castelnuovo polynomial. 
\end{proof} 
\begin{proof}[Proof of Corollary \ref{ref-1.6-10}] 
The number of solutions to the conditions (1-3) in the statement of 
Theorem C is 
\[ 
\prod_{l>\sigma} \left(\biggl(\sum_{i\le 
    l}q_i\biggr)-\max(q_l,0)\right)= 
\prod_{l>\sigma}\min\biggl( \sum_{i< l}q_i, 
\sum_{i\le l}q_i\biggr) 
\] 
Noting that $\sum_{i\le l}q_i=1+s_{l-1}-s_{l}$ finishes the proof. 
\end{proof} 

\begin{convention} 
 Below we will call a formal power series of the form 
\[ 
\frac{1}{(1-t)^3}-\frac{s(t)}{1-t} 
\] 
where $s(t)$ is a Castelnuovo polynomial of weight $n$ an 
\emph{admissible Hilbert series of weight $n$}. 
\end{convention} 
\section{The stratification by Hilbert series} 
In this section we will prove the following result. 
\begin{theorem} \label{ref-7.1-42} There is a (weak) stratification into 
smooth, non-empty connected 
  locally closed sets 
\begin{equation} 
\label{ref-7.1-43} 
\Hilb_n(\PP^2_q)=\bigcup_h \Hilb_h(\PP^2_q) 
\end{equation} 
where the union runs over the (finite set) of admissible Hilbert 
series of weight $n$ and where the points in $\Hilb_h(\PP^2_q)$ 
represents the points in $\Hilb_n(\PP^2_q)$ corresponding to objects with 
Hilbert series $h$. 

Furthermore we have 
\begin{equation} 
\label{ref-7.2-44} 
\overline{\Hilb_h(\PP^2_q)}\subset \bigcup_{h'\ge h} \Hilb_{h'}(\PP^2_q) 
\end{equation} 
In the decomposition \eqref{ref-7.1-43} there is a unique stratum of 
maximal dimension $2n$ which corresponds to the Hilbert series 
$h_{\text{min}}(t)$  (see \S\ref{ref-1-0}). 
\end{theorem} 
That the strata are non-empty is Theorem A. The rest of 
Theorem \ref{ref-7.1-42} will be a consequence of 
Lemma \ref{ref-7.1.1-49}, 
Corollary \ref{ref-7.2.3-58} 
and Proposition \ref{ref-7.3.1-60} 
below. 
\medskip 

We refer to \eqref{ref-7.1-43} as a ``weak'' stratification (for an 
ordinary stratification one would require the inclusions in 
\eqref{ref-7.2-44} to be equalities, which is generally not the 
case). 
\medskip 

 In the commutative case Theorem \ref{ref-7.2-44} was proved by 
Gotzmann \cite{Gotzmann}. It is not clear to us that Gotzmann's method can be 
generalized to the non-commutative case. In any case, the reader will 
notice, that our proof is substantially different. 
\begin{proof}[Proof of Theorem B] 
This is now clear from Theorem \ref{ref-7.1-42}. 
\end{proof} 
It follows from Theorem C that given a Hilbert 
series $h(t)=q(t)/(1-t)^3$ there is a unique legal choice of Betti numbers 
$(a_i)_i$, $(b_i)_i$ such that $a_i$ and $b_i$ are not both non-zero 
for all $i$. Namely 
\begin{equation} 
\label{ref-7.3-45} 
(a_i,b_i)= 
\begin{cases} 
(q_i,0)&\text{if $q_i\ge 0$}\\ 
(0,-q_i)&\text{otherwise} 
\end{cases} 
\end{equation} 
We call this the \emph{minimal Betti numbers} associated to $h$. 

We have some extra information on the strata $\Hilb_h(\PP^2_q)$. Define 
$\Hilb_h(\PP^2_q)^{\min}$ as the subset of $\Hilb_h(\PP^2_q)$ 
consisting of objects with minimal Betti numbers. 
\begin{proposition} 
\label{ref-7.2-46} 
$\Hilb_h(\PP^2_q)^{\min}$ is open in $\Hilb_h(\PP^2_q)$ 
\end{proposition} 
This is proved in \S\ref{ref-7.1-48} below. 

\medskip 

Assume that $A$ is elliptic and that in the geometric data 
$(E,\Oscr_E(1),\sigma)$ associated to $A$, $\sigma$ has infinite 
order.  Let $\Hilb_n(\PP^2_q)^{\operatorname{inv}}$ be the reflexive 
objects in $\Hilb_n(\PP^2_q)$. This is an open subset (see 
\cite[Theorem 8.11]{NS}). 
\begin{proposition} \label{ref-7.3-47}  For all admissible Hilbert series 
$h$ with 
  weight $n$ we have 
\[ 
\Hilb_h(\PP^2_q)\cap \Hilb_n(\PP^2_q)^{\operatorname{inv}} \neq 
\emptyset 
\] 
\end{proposition} 
\begin{proof} 
  This is a special case of Proposition \ref{ref-5.4.1-41}. 
\end{proof} 
\begin{remark} 
  Consider the Hilbert scheme of points $\Hilb_{n}(\PP^{2})$ in the 
  projective plane $\PP^{2}$.  The inclusion relation between the 
  closures of the strata of $\Hilb_{n}(\PP^{2})$ has been a subject of 
  interest in \cite{BH, maC, CW, HRW}.  Although in general the 
  precise inclusion relation is still unknown, the special case where 
  the Hilbert series of the strata are as close as possible is 
  completely settled (see \cite{DV2, Guerimand}).  It is a natural to 
  consider the same question for the varieties $\Hilb_{n}(\PP_q^{2})$, 
  where one may use the same techniques as in \cite{DV2}. 
\end{remark} 
\subsection{Moduli spaces} 
\label{ref-7.1-48} 
In this section ``points'' of schemes will be not necessarily closed. 
We will consider functors from the category of 
noetherian $k$-algebras $\Noeth/k$ to the category of sets. For $R\in 
\Noeth/k$ we write $(-)_R$ for the base extension $-\otimes R$. If 
$x$ is a (not necessarily closed) point in $\Spec R$ then we write 
$(-)_x$ for the base extension $-\otimes_R k(x)$. We put 
$\PP^2_{q,R}=\Proj A_R$. 

It follows from \cite[Prop.\ 4.9(1) and 4.13]{ASZ} that $A$ is 
strongly noetherian so $A_R$ is still noetherian. Furthermore it 
follows from \cite[Prop.\ C6]{AZ1} that $A_R$ satisfies the 
$\chi$-condition and finally by \cite[Cor.\ C7]{AZ1} 
$\Gamma(\PP^2_{q,R},-)$ has cohomological dimension two.

An $R$-family of objects in $\coh(\PP^2_q)$ or $\grmod(A)$ is by 
definition an $R$-flat object \cite{AZ1} in these categories. 

For $n\in \NN$ let 
$\Hscr\mathit{ilb}_n(\PP^2_q)(R)$ be the $R$-families of objects 
$\Iscr$ in $\coh(\PP^2_q)$, \emph{modulo Zariski local isomorphism on $\Spec 
R$}, with the property that for any map $x\in \Spec R$, $\Iscr_x$ is 
torsion free normalized of rank one in $\coh(\PP^2_{q,k(x)})$. 

The main result of \cite{NS} is that $\Hscr\mathit{ilb}_n(\PP^2_q)$ is 
represented by a smooth scheme $\Hilb_{n}(\PP^2_q)$ of dimension $2n$ 
(see also \cite{DV} for a special case, treated with a different 
method which yields some extra information). 
\begin{warning} The reader will notice that now the set 
  $\Hilb_{n}(\PP^2_q)(k)=\Hscr\mathit{ilb}_n(\PP^2_q)(k)$ parametrizes 
objects in $\coh(\PP^2_q)$ rather than in $\grmod(A)$ as was the case 
in the introduction. However by Corollary \ref{ref-3.4.2-22} the new 
point of view is equivalent to the old one. 
\end{warning} 
If $h(t)$ is a admissible Hilbert series of weight $n$ then 
$\Hscr\mathit{ilb}_{h}(\PP^2_q)(R)$ is the set of $R$-families of 
torsion free graded $A$-modules which have Hilbert series $h$ and 
which have projective dimension one, \emph{modulo local isomorphism on 
  $\Spec R$}. The map $\pi$ defines a map 
\[ 
\pi(R):\Hscr\mathit{ilb}_{h}(\PP^2_q)(R)\r 
\Hscr\mathit{ilb}_n(\PP^2_q)(R): I\mapsto \pi I 
\] 
 Below we will write 
$\Iscr^u$ for a universal family on $\Hilb_{n}(\PP^2_q)$. This is a 
sheaf of graded $\Oscr_{\Hilb_{n}(\PP^2_q)}\otimes A$-modules on 
$\Hilb_{n}(\PP^2_q)$. 
\begin{lemmas} 
\label{ref-7.1.1-49} 
  The map $\pi(k)$ is an injection which identifies 
  $\Hscr\mathit{ilb}_{h}(\PP^2_q)(k)$ with 
\[ 
\{ 
x\in \Hscr\mathit{ilb}_{n}(\PP^2_q)(k)\mid h_{\Iscr^u_x}=h 
\} 
\] 
This is a  locally closed subset of 
$\Hscr\mathit{ilb}_{n}(\PP^2_q)(k)$. Furthermore 
\begin{equation} 
\label{ref-7.4-50} 
\overline{\Hscr\mathit{ilb}_h(\PP^2_q)(k)}\subset \bigcup_{h'\ge h} 
\Hscr\mathit{ilb}_{h'}(\PP^2_q)(k) 
\end{equation} 
\end{lemmas} 
\begin{proof}  The fact that $\pi(k)$ is an injection 
and does the required identification follows from 
  Corollary \ref{ref-3.4.2-22}. 

For any $N\ge 0$ we have by 
Corollary \ref{ref-A.3-64} that 
\[ 
\Hscr\mathit{ilb}_{h,N}(\PP^2_q)(k)=\{x\in 
\Hscr\mathit{ilb}_{n}(\PP^2_q)(k)\mid h_{\Iscr^u_x}(n)= 
h(n)\text{ for }n\le N\} 
\] 
is locally closed in $\Hscr\mathit{ilb}_{n}(\PP^2_q)(k)$. By 
Theorem A we know that only a finite 
number of Hilbert series occur for objects in 
$\Hscr\mathit{ilb}_{n}(\PP^2_q)(k)$. Thus 
$\Hscr\mathit{ilb}_{h,N}(\PP^2_q)(k)=\Hscr\mathit{ilb}_{h}(\PP^2_q)(k)$ 
for $N\gg 0$. \eqref{ref-7.4-50} also follows easily from semi-continuity. 
\end{proof} 
Now let $\Hilb_{h}(\PP^2_q)$ be the reduced locally closed subscheme 
of $\Hilb_{n}(\PP^2_q)$ whose closed points are given by 
$\Hscr\mathit{ilb}_{h}(\PP^2_q)(k)$. We then have the following 
result. 
\begin{propositions} \label{ref-7.1.2-51} $\Hilb_{h}(\PP^2_q)$ represents 
the functor 
  $\Hscr\mathit{ilb}_{h}(\PP^2_q)$. 
\end{propositions} 
Before proving this proposition we need some technical results. 
The following is proved in \cite{NS}. For the 
convenience of the reader we put the proof here. 
\begin{lemmas} Assume that 
  $\Iscr$, $\Jscr$ are $R$-families of objects in 
  $\coh(\PP^2_q)$ with the property that for any map $x\in \Spec R$, 
  $\Iscr_x$ is torsion free of rank one in 
  $\coh(\PP^2_{q,k(x)})$.  Then $\Iscr$, $\Jscr$ represent the same 
  object in $\Hscr\mathit{ilb}_n(\PP^2_q)(R)$ 
if and only if there is an invertible module $\frak{l}$ in 
  $\Mod(R)$ such that 
\[ 
\Jscr=\frak{l}\otimes_R \Iscr 
\] 
\end{lemmas} 
\begin{proof} Let $\Iscr$ be as in the statement of the lemma. We 
  first claim that the natural map 
\begin{equation} 
\label{ref-7.5-52} 
R\r \End(\Iscr) 
\end{equation} 
is an isomorphism. Assume first that $f\neq 0$ is in the kernel of 
\eqref{ref-7.5-52}. Then the flatness of $\Iscr$ implies $\Iscr \otimes_R 
Rf=0$. This implies that $\Iscr_x=\Iscr\otimes_R k(x)=0$ for some 
$x\in \Spec R$ and this is a contradiction since by definition 
$\Iscr_x\neq 0$. 

It is easy to see \eqref{ref-7.5-52} is surjective (in fact an isomorphism) 
when $R$ is a field. It follows that for all $x\in \Spec R$ 
\[ 
\End(\Iscr)\otimes_R k(x)\r \End(\Iscr\otimes_R k(x)) 
\] 
is surjective. Then it follows from base change (see \cite[Thm 
4.3(1)(4)]{NS}) that $\End(\Iscr)\otimes_R k(x)$ is one dimensional 
and hence 
\eqref{ref-7.5-52} is surjective by Nakayama's lemma. 

Now let $\Iscr$, $\Jscr$ be as in the statement of the lemma and 
assume they represent the same element of 
$\Hscr\mathit{ilb}_n(\PP^2_q)(R)$, i.e.\ they are locally isomorphic. Put 
\[ 
\frak{l}=\Hom(\Iscr,\Jscr) 
\] 
It is easy to see that $\frak{l}$ has the required properties since 
this may be checked locally on $\Spec R$ and then we may invoke the 
isomorphism \ref{ref-7.5-52}. 
\end{proof} 
\begin{lemmas} 
\label{ref-7.1.4-53} 
  Assume that $R$ is finitely generated and let $P_0$, $P_1$ be 
  finitely generated graded free $A_R$-modules. Let $N\in 
  \Hom_A(P_1,P_0)$. Then 
\[ 
V=\{x\in \Spec R\mid N_x\text{ is injective  with torsion free cokernel}\} 
\] 
is open. Furthermore the restriction of $\coker N$ to $V$ is $R$-flat. 
\end{lemmas} 
\begin{proof} 
We first note that the formation of $V$ is compatible with base 
change. It is sufficient to prove this for an extension of fields. The 
key point is that if $K\subset L$ is an extension of fields and $M\in 
\grmod(A_K)$ then $M$ is torsion free if and only if $M_L$ is torsion 
free. This follows from the fact that if $D$ is the graded quotient 
field of $A_K$ then $M$ is torsion free if and only if the map $M\r 
M\otimes_{A_K} D$ is injective.

To prove openness of $V$ we may now assume by \cite[Theorem 0.5]{ASZ} 
that $R$ is a Dedekind domain (not necessarily finitely generated).

Assume $K=\ker N\neq 0$. Since $\gldim R=1$ we deduce that the map $K\r 
P_1$ is  degree wise split. Hence $N_x$ is never injective and 
the set $V$ is empty. 

So we assume $K=0$ and we let $C=\coker N$. Let $T_0$ be the 
$R$-torsion part of $C$. Since $A_R$ is noetherian $T_0$ is finitely 
generated. We may decompose $T_0$ degree wise according to the 
maximal ideals of $R$. Since it is clear that this yields a decomposition of 
$T_0$ as $A_R$-module it follows that  there can be only a finite number of 
points in the support of $T_0$ as $R$-module. 

If $x\in \Spec R$ is in the support of $T_0$ as $R$-module then 
$\Tor_1^R(C,k(x))\neq 0$ and hence $N_x$ is not injective. Therefore 
$x\notin V$.  By considering an affine covering of the complement of 
the support of $T_0$ as $R$-module we reduce to the case where $C$ is 
torsion free as $R$-module. 

Let $\eta$ be the generic point of $\Spec R$ and assume that $C_\eta$ 
has a non-zero torsion submodule $T_\eta$. Put $T=T_\eta\cap C$. Since 
$R$ is Dedekind the map $T\r C$ is degree wise split.  Hence 
$T_{k(x)}\subset C_{k(x)}$ and so $C_{k(x)}$ will always have 
torsion. Thus $V$ is empty. 

Assume $T_\eta=0$. It is now sufficient to construct an non-empty open 
$U$ in $\Spec R$ such that $U\subset V$.  We have an embedding 
$C\subset C^{\ast\ast}$.  Let $Q$ be the maximal $A_R$ submodule of 
$C^{\ast\ast}$ containing $C$ such that such that $Q/C$ is 
$R$-torsion. Since $Q/C$ is finitely generated it is supported on a 
finite number of closed points of $\Spec R$ and we can get rid of 
those by considering an affine open of the complement of those 
points.

Thus we may assume that $C^{\ast\ast}/C$ is $R$-torsion free. Under 
this hypothesis we will prove that $C_x$ is torsion free for all 
closed points $x\in \Spec R$. Since we now have an injection $C_x\r 
(C^{\ast\ast})_x$ it is sufficient to prove that $(C^{\ast\ast})_x$ is 
torsion free. To this end we way assume that $R$ is a discrete 
valuation ring and $x$ is the closed point of $\Spec R$. 

Let $\Pi$ be the uniformizing element of $R$ and let $T_1$ be the 
torsion submodule of $(C^{\ast\ast})_x$. Assume $T_1\neq 0$ and let 
$Q$ be its inverse image in $C^{\ast\ast}$. Thus we have an exact 
sequence 
\begin{equation} 
\label{ref-7.6-54} 
0\r \Pi C^{\ast\ast}\r Q\r T_1\r 0 
\end{equation} 
which is cannot be split since otherwise $T_1\subset C^{\ast\ast}$ 
which is impossible. 

We now apply $(-)^\ast$ to \eqref{ref-7.6-54}. Using 
$\underline{\Ext}^1_{A_R}(T_1, A_R)=\underline{\Hom}_{A_x}(T_1, 
A_x)=0$ we deduce $Q^\ast=C^{\ast\ast\ast}=C^\ast$. Applying 
$(-)^\ast$ again we deduce $Q^{\ast\ast}=C^{\ast\ast}$ and hence the 
map $Q\r Q^{\ast\ast}\cong C^{\ast\ast}$ gives a splitting of 
\eqref{ref-7.6-54}, which is a contradiction. This finishes the proof of 
the openness of $V$. 

The flatness assertion may be checked locally. So we may assume that 
$R$ is a local ring with closed point $x$ and $x\in V$.  Thus for any 
$m$ we have a map between free $R$-modules $(P_1)_m\r (P_0)_m$ which 
remains injective when tensored with $k(x)$. A standard application 
of Nakayama's lemma then yields that the map is split, and hence its 
cokernel is projective. 
\end{proof} 
\begin{lemmas} 
\label{ref-7.1.5-55} Assume $I\in \Hscr\mathit{ilb}_{h}(\PP^2_q)(R)$ 
  and $x\in \Spec R$. Then there exist: 
\begin{enumerate} 
\item an element $r\in R$ with $r(x)\neq 0$; 
\item a polynomial ring 
  $S=k[x_1,\ldots, x_n]$; 
\item a point $y\in \Spec S$; 
\item an element $s\in S$ with $s(y)\neq 0$; 
\item a homomorphism of rings $\phi:S_s\r R_r$ such that $\phi(x)=y$ 
(where we also have written $\phi$ for the dual map $\Spec R_r\r \Spec 
S_s$); 
\item an object $I^{(0)}$ in  $\Hscr\mathit{ilb}_{h}(\PP^2_q)(S_s)$ 
  such that $I^{(0)}\otimes_{S_s} R_r=I\otimes_{S} S_s$. 
\end{enumerate} 
\end{lemmas} 
\begin{proof} 
By hypotheses $I$ has a presentation 
\[ 
0\r P_1\r P_0\r I\r 0 
\] 
where $P_0$, $P_1$ are finitely graded projective $A_R$-modules. It is 
classical that we have $P_0\cong \frak{p}_0 \otimes_R A$, $P_1\cong 
\frak{p}_1 \otimes_R A$ where $\frak{p}_0$, $\frak{p}_1$ are finitely 
generated graded projective $R$-modules.  By localizing $R$ at an 
element which is non-zero in $x$ we may assume that $P_0$, $P_1$ are 
graded free $A_R$-modules. After doing this $N$ is given by a $p\times 
q$-matrix with coefficients in $A_R$ for certain $p$, $q$. 

 Then by choosing a 
  $k$-basis for $A$ and writing out the entries of $N$ in terms of this basis 
  with coefficients in $R$ we may construct a polynomial ring 
 $S=k[x_1,\ldots, x_n]$ together with a morphism $S\r R$ and a 
 $p\times q$-matrix $N^{(0)}$ over $A_S$ such that $N$ is obtained by 
  base-extension from $N^{(0)}$. Thus $I$ is obtained by 
  base-extension from the cokernel $I^{(0)}$ of a map 
\[ 
N^{(0)}:P_1^{(0)}\longrightarrow P_0^{(0)} 
\] 
where $P_1^{(0)}$, $P_0^{(0)}$ are graded free $A_S$-modules.  Let 
$y$ be the image of $x$ in $\Spec S$. By construction we have $I_x= 
I^{(0)}_y\otimes_{k(y)} k(x)$. From this it easily follows that 
$I^{(0)}_y\in \Hscr\mathit{ilb}_h(k(y))$. 

The module 
$I^{(0)}$ will not in general satisfy the requirements of the 
lemma but it follows from Lemma \ref{ref-7.1.4-53} that this will be the case 
after inverting a suitable element in $S$ non-zero in $y$. This 
finishes the proof. 
\end{proof} 
\begin{proof}[Proof of Proposition \ref{ref-7.1.2-51}] 
  Let $R\in \Noeth/k$. We will construct inverse bijections 
\begin{gather*} 
\Phi(R):\Hscr\mathit{ilb}_{h}(\PP^2_q)(R) \r \Hom(\Spec 
R,\Hilb_{h}(\PP^2_q))\\ 
\Psi(R):\Hom(\Spec 
R,\Hilb_{h}(\PP^2_q))\r \Hscr\mathit{ilb}_{h}(\PP^2_q)(R) 
\end{gather*} 
We start with $\Psi$. 
For  $w\in \Hom(\Spec 
R,\Hilb_{h}(\PP^2_q))$  we put 
\[ 
 \Psi(R)(w)=\omega(\Iscr^u_R)=\bigoplus_m \Gamma(\PP^2_{q,R},\Iscr^u_R(m)) 
\] 
We need to show that $\omega(\Iscr^u_R)\in 
\Hscr\mathit{ilb}_{h}(\PP^2_q)(R)$. It is 
clear that this can be done Zariski locally on $\Spec R$. Therefore we 
may assume that $w$ factors as 
\[ 
\Spec R \r \Spec S\r \Hilb_{h}(\PP^2_q) 
\] 
where $\Spec S$ is an affine open subset of $\Hilb_{h}(\PP^2_q)$. 

Now by lemma \ref{ref-7.1.1-49} 
\[ 
\Spec S\r \NN: x\mapsto \dim_k \Gamma(\PP^2_{q,x},\Iscr^u_x(m)) 
\] 
has constant value $h(m)$ and hence by Corollary \ref{ref-A.4-66} below 
$\Gamma(\PP^2_{q,S},\Iscr^u_S(m))$ is a projective  $S$-module and 
furthermore by \cite[Lemma C6.6]{AZ1}. 
\begin{align*} 
  \Gamma(\PP^2_{q,x},\Iscr^u_x(m))&=\Gamma(\PP^2_{q,S},\Iscr^u_S(m)) 
  \otimes_S k(x)\\ 
  \Gamma(\PP^2_{q,R},\Iscr^u(m)_R)&=\Gamma(\PP^2_{q,S},\Iscr^u_S(m))\otimes_S 
  R 
\end{align*} 
for $x\in \Spec S$. We deduce that $\omega(\Iscr^u_S)$ is flat and furthermore 
\begin{align*} 
\omega(\Iscr^u_S)_x&=\omega(\Iscr^u_x)\\ 
\omega(\Iscr^u_S)_R&=\omega(\Iscr^u_R) 
\end{align*} 
Using the first equation we deduce from Corollary \ref{ref-3.4.2-22} 
and Nakayama's lemma that $\omega(\Iscr^u_S)$ has projective dimension 
one. Thus  $\omega(\Iscr^u_S)\in 
\Hscr\mathit{ilb}_{h}(\PP^2_q)(S)$. From the second equation we then 
deduce $\omega(\Iscr^u_R)\in 
\Hscr\mathit{ilb}_{h}(\PP^2_q)(R)$. 

Now we define $\Phi$. Let $I\in \Hscr\mathit{ilb}_h(R)$.  We define 
$\Phi(R)(I)$ as the map $w:\Spec R\r \Hilb_n(\PP^2_q)$ corresponding 
to $\Pi I$. I.e.\ formally 
\[ 
\pi I=\Iscr^u_{w}\otimes_R \frak{l} 
\] 
where $\frak{l}$ is an invertible $R$ module and where this time we 
have made the base change map $w$ explicit in the notation.  We need 
to show that $\im w$ lies in $\Hilb_h(\PP^2_q)$. Again we may do this 
locally on $\Spec R$. Thus by lemma \ref{ref-7.1.5-55} we may assume that 
there is a map $\theta:S\r R $ where $S$ is integral and finitely 
generated over $k$ and $I$ is obtained from $I^{(0)}\in 
\Hscr\mathit{ilb}_h(S)$ by base change. Let $v:\Spec S\r 
\Hilb_n(\PP^2_q)$ be the map corresponding to $I^{(0)}$. An elementary 
computation shows that  $v\theta=w$. In other words it is sufficient 
to check that $\im v\subset \Hilb_h(\PP^2_q)$. But since $S$ is 
integral of finite type over $k$ it suffices to check this for 
$k$-points. But then it follows from Lemma \ref{ref-7.1.1-49}. 

We leave to the reader the purely formal computation that $\Phi$ and 
$\Psi$ are each others inverse. 
\end{proof} 
\begin{proof}[Proof of Proposition \ref{ref-7.2-46}] Let $(a_i)_i$, 
  $(b_i)_i$ be minimal Betti numbers corresponding to $h$. Let $I^u$ be 
  the universal family on $\Hilb_h(\PP^2_q)$. Then it is easy to see 
  that 
\[ 
\Hilb_h(\PP^2_q)^{\min}=\{x\in \Hilb_h(\PP^2_q) 
\mid \forall i: \dim_{k(x)} (I_x^u\otimes_{A_x} k(x))_i=a_i\} 
\] 
It follows from Lemma \ref{ref-A.1-62} below that this defines an open 
subset. 
\end{proof} 
\subsection{Dimensions} 
Below a point will again be be closed point. 
\begin{lemmas} Let $I\in \Hilb_h(\PP^2_q)$. Then 
  canonically 
\[ 
T_I(\Hilb_h(\PP^2_q))\cong \Ext^1_A(I,I) 
\] 
\end{lemmas} 
\begin{proof} 
If $\Fscr$ is a functor from (certain) rings to sets and $x\in 
\Fscr(k)$ then the tangent space $T_x(\Fscr)$ is by definition the inverse 
image of $x$ under  the map 
\[ 
\Fscr(k[\epsilon]/(\epsilon^2))\r \Fscr(k) 
\] 
which as usual is canonically a $k$-vector space.  If 
$\Fscr$ is represented by a scheme $F$ then of course 
$T_x(\Fscr)=T_x(F)$. 

The proposition follows from the fact that if  $I\in 
\Hscr\mathit{ilb}_h(\PP^2_q)(k)$ then the 
tangent space $T_I(\Hscr\mathit{ilb}_h(\PP^2_q))$ is canonically 
identified with $\Ext^1_A(I,I)$ (see \cite[Prop.\ E1.1]{AZ1}). 
\end{proof} 
We  now express 
  $\dim_{k}\Ext^1_A(I,I)$ in terms of 
$s_{I}(t)$. 
\begin{propositions} 
\label{ref-7.2.2-56} 
Let $I\in \Hilb_h(\PP^2)$ and assume $I\neq A$. Let $s_{I}(t)$ be the 
Castelnuovo polynomial of~$I$. Then we have 
\[ 
\dim_k \Ext^1_A(I,I)=1+n+c 
\] 
where $n$ is the invariant of $I$ and $c$ is the constant term of 
\begin{equation} 
\label{ref-7.7-57} 
(t^{-1}-t^{-2})s_I(t^{-1})s_I(t) 
\end{equation} 
In particular this dimension is independent of $I$. 
\end{propositions} 
\begin{corollarys} \label{ref-7.2.3-58} $\Hilb_h(\PP^2_q)$ is smooth of 
dimension $1+n+c$ 
  where $c$ is as in the previous theorem. 
\end{corollarys} 
\begin{proof} This follows from the fact that the tangent spaces of 
  $\Hilb_h(\PP^2_q)$  have constant dimension $1+n+c$. 
\end{proof} 
\begin{proof}[Proof of Proposition  \ref{ref-7.2.2-56}] 
We start with the following observation. 
\[ 
\sum_i 
(-1)^ih_{\underline{\Ext}^i_A(M,N)}(t)=h_{M}(t^{-1})h_{N}(t)(1-t^{-1})^{3} 
\] 
for $M,N \in \grmod(A)$. 
This follows from the fact that both sides a additive on short exact 
sequences, and they are equal for $M=A(-i)$, $N=A(-j)$. 

Applying this with $M=N=I$ and using the fact that $\pdim I=1$, 
$\underline{\Hom}_A(I,I)=k$ 
we obtain that $\dim_k \Ext^1_A(I,I)$ is the constant term of 
\begin{equation*} 
\begin{split} 
1-h(t^{-1})h(t)(1-t^{-1})^3 
& = 1-(1-t^{-1})^3\left(\frac{1}{(1-t^{-1})^3}-\frac{s(t^{-1})} 
{1-t^{-1}}\right)\left(\frac{1}{(1-t)^3}-\frac{s(t)}{1-t}\right) \\ 
& = 1-\frac{1}{(1-t)^3}+\frac{s(t)}{1-t}+\frac{t^{-2}s(t^{-1})}{1-t} 
-t^{-2}(1-t)s(t^{-1})s(t) 
\end{split} 
\end{equation*} 
(where we dropped the index ``$I$''). 
Introducing the known constant terms finishes the proof. 
\end{proof} 
\begin{corollarys} 
Let $\Iscr \in \Hilb_n(\PP^2_q)$ and $I=\omega\Iscr$. 
Then 
\begin{equation} 
\label{ref-7.8-59} 
\min(n+2,2n) \leq \dim_k \Ext^1_A(I,I) \le 2n 
\end{equation} 
with  equality on the left if and only if $h_I(t)=h_{\max}(t)$ and 
equality on the right if and only if $h_I(t)=h_{\min}(t)$ (see \S\ref{ref-1-0}). 
\end{corollarys} 
\begin{proof}  Since the case $n=0$ is obvious we assume below $n\ge 1$. 
We compute the constant term of \eqref{ref-7.7-57}. 
Put $s(t)=s_I(t)=\sum s_i t^i$. 
Thus the sought constant term is the difference between the coefficient 
of $t$ and the coefficient of $t^2$ in $s(t^{-1})s(t)$. 
This difference is 
\[ 
\sum_{j-i=1} s_i s_j-\sum_{j-i=2}s_is_j 
\] 
which may be rewritten as 
\begin{align*} 
\sum_{j} s_{j+1} s_j-\sum_{j}s_{j+2}s_j & 
=\sum_{j} s_{j} s_{j-1}-\sum_{j}s_{j}s_{j-2} 
=\sum_j s_j(s_{j-1}-s_{j-2}) 
\end{align*} 
Now we always have $s_{j-1}-s_{j-2}\le 1$ and $s_{-1}-s_{-2}=0$. Thus 
\[ 
\sum_j s_j(s_{j-1}-s_{j-2})\le -1+\sum_j s_j=n-1 
\] 
which implies $\dim_k \Ext^1_A(I,I) \le 2n $ by Proposition 
\ref{ref-7.2.2-56}, and 
we will clearly have equality if and only if $s_{j-1}-s_{j-2}=1$ for 
$j>0$ and $s_j\neq 0$. This is equivalent to $s(t)$ being of the 
form 
\[ 
1+2t+3t^2+\cdots+(u-1)t^u+vt^{u+1} 
\] 
for some integers $u > 0$ and $v\ge 0$. 
This in turn is equivalent with $h_I(t)$ being equal to $h_{\min}(t)$. 
This proves the upper bound of $\eqref{ref-7.8-59}$. \\ 

Now we prove the lower bound. Since $s_I(t)$ is a Castelnuovo 
polynomial it has the form 
\[ 
s(t) = 1 + 2t + 3t^{2} + \ldots + \sigma t^{\sigma-1} + 
s_{\sigma}t^{\sigma} + s_{\sigma + 1}t^{\sigma + 1} + \ldots 
\] 
where 
\[ 
\sigma \geq s_{\sigma} \geq s_{\sigma + 1} \geq \ldots 
\] 
We obtain 
\begin{align*} 
c & = \sum_{j} s_{j}(s_{j+1} - s_{j+2}) \\ 
& = - (1 +  2 + 3 + \ldots + (\sigma - 2))  + 
\sum_{j \geq \sigma - 2}s_{j}(s_{j+1} - s_{j+2}) 
\end{align*} 
We denote the subsequence obtained by dropping the zeros from the sequence 
of non-negative integers 
$\left( s_{j+1} - s_{j+2} \right)_{j \geq \sigma -2}$ by $e_{1}, e_{2}, 
\ldots, e_{r}$. 
Note that $\sum_{i}e_{i} = \sigma$. 
We get 
\begin{align*} 
c \geq & - (1 + 2 + 3 + \ldots + (\sigma - 2)) \\ 
& + (\sigma - \delta) e_{1} + (\sigma - e_{1})e_{2} + \ldots + 
(\sigma - e_{1} - \ldots - e_{r-1})e_{r} 
\end{align*} 
where $\delta = 1$ if $s_{\sigma} < \sigma$ and $0$ otherwise. 
Now we have that 
\begin{gather*} 
(\sigma - e_{1} - \cdots - e_{r-1})e_{r} = e_{r} e_{r}  \geq 1 + \cdots + 
e_{r} \\ 
(\sigma - e_{1} - \cdots - e_{r-2})e_{r-1} = (e_{r-1} + e_{r}) e_{r-1} 
\geq (1+e_{r}) + \cdots + (e_{r-1} + e_{r}) \\ 
\vdots \\ 
(\sigma - e_{1})e_{2} = (e_{2} + \cdots + e_{r}) e_{2}  \geq 
(1 + e_{3} + \cdots + e_{r}) + \cdots + (e_{2} + \cdots + e_{r}) \\ 
\sigma e_{1} = (e_{1} + \cdots + e_{r}) e_{1}  \geq 
(1 + e_{2} + \cdots + e_{r}) + \cdots + (e_{1} + \cdots + e_{r}) 
\end{gather*}

hence 
\begin{align*} 
c \geq 2\sigma - 1 - \delta e_{1} 
\end{align*} 
Hence $c \geq 0$ and $c = 0$ if and only if $\sigma = 1$, $r = 1$ and 
$\delta = 1$, so if and only if $s_{I}(t) = 1$. In that case, the invariant 
$n$ of $I$ is $1$. 
If $n > 1$ then $c \geq 1$ which proves the lower bound of 
\eqref{ref-7.8-59} by Proposition \ref{ref-7.2.2-56}. 
Clearly $c = 1$ if and only if $\sigma = 1$ and $r=1$, which 
is equivalent with $h_I(t)$ being equal to $h_{\max}(t)$. 
\end{proof} 
\begin{remarks} 
  The fact that $\dim_k \Ext^1_A(I,I)\le 2n$ can be shown 
  directly. Indeed from the formula 
\[ 
\Ext^1_{\Tails(A)}(\Iscr,\Iscr) \cong \lim_{\lr} 
\Ext^1_A(I_{\ge n},I) 
\] 
and from the fact that $\Ext^1_A(k,I)=0$ we obtain an injection 
\[ 
\Ext^1_A(I,I)\hookrightarrow \Ext^1_{\Tails(A)}(\Iscr,\Iscr) 
\] 
and the right hand side is the tangent space $\Iscr$ in the smooth 
variety $\Hilb_n(\PP^2_q)$ which has dimension $2n$. 
\end{remarks} 
\subsection{Connectedness} 
In this section we prove 
\begin{propositions} 
\label{ref-7.3.1-60} Assume that $h$ is an admissible Hilbert 
  polynomial. Then any two points in $\Hilb_h(\PP^2_q)$ can be 
  connected using an open subset of an affine line. 
\end{propositions} 
\begin{proof} 
Let $I,J\in \Hilb_h(\PP^2_q)$. Then $I,J$ have resolutions 
\begin{gather*} 
0\r \oplus_i A(-i)^{b_i}\r \oplus_i A(-i)^{a_i} \r I\r 0\\ 
0\r \oplus_i A(-i)^{d_i}\r \oplus_i A(-i)^{c_i} \r J\r 0 
\end{gather*} 
where $a_i-b_i=c_i-d_i$. Adding terms of the form 
$A(-j)\xrightarrow{\Id} A(-j)$ we may change these resolutions to have 
the following form 
\begin{gather*} 
0\r \oplus_i A(-i)^{f_i}\xrightarrow{M} \oplus_i A(-i)^{e_i}\r I\r 0\\ 
0\r \oplus_i A(-i)^{f_i}\xrightarrow{N} \oplus_i A(-i)^{e_i}\r J\r 0 
\end{gather*} 
for matrices $M,N\in H=\Hom_A(\oplus_i A(-i)^{e_i},\oplus_i 
A(-i)^{f_i})$. Let $L\subset H$ be the line through $M$ and $N$. Then 
by Lemma \ref{ref-7.1.4-53} an open set of $L$ defines points in 
$\Hilb_h(\PP^2_q)$. This finishes the proof. 
\end{proof} 

\appendix 
\section{Upper semi-continuity for non-commutative $\Proj$} 
\label{ref-A-61} 
In this section we discuss some results which are definitely at least 
implicit in 
\cite{AZ1} but for which the authors have been unable to find a 
convenient reference.  The methods are quite routine. 
We refer to \cite{AZ1,H} for more details. 

Below $R$ will be a noetherian commutative ring and 
$A=R+A_1+A_2+\cdots$ is a noetherian connected graded $R$-algebra. 
\begin{lemma} 
\label{ref-A.1-62} 
Let $M\in \grmod(A)$ be flat over $R$ and $n\in \ZZ$. Then the function 
\[ 
\Spec R\r \ZZ:x\mapsto \dim_k \Tor_i^{A_{k(x)}}(M_{k(x)},k(x))_n 
\] 
is upper semi-continuous. 
\end{lemma} 
\begin{proof} Because of flatness we have 
$\Tor_i^{A_{k(x)}}(M_{k(x)},k(x))=\Tor_i^{A}(M,k(x))$. 
Let $F^\cdot\r M\r 0$ be a graded  resolution of $M$ consisting of 
free $A$-modules of finite rank. Then 
$\Tor_i^{A}(M,k(x))_n$ is the homology of $(F^\cdot)_n\otimes_A 
{k(x)}$. Since $(F^\cdot)_n$ is a complex of free $R$-modules, the 
result follows in the usual way. 
\end{proof}

Now  we 
write $X=\Proj A$ and we use the associated notations as outlined in 
\S\ref{ref-3.1-15}.  In addition we will assume that $A$ satisfies the 
following conditions. 
\begin{enumerate} 
\item $A$ satisfies $\chi$ \cite{AZ}. 
\item $\Gamma(X,-)$ has finite cohomological dimension. 
\end{enumerate} 
Under these hypotheses we prove 
\begin{proposition} 
\label{ref-A.2-63} 
Let $\Gscr\in \coh(X)$ be flat over $R$ and let 
  $\Fscr\in \coh(X)$ be arbitrary. Then there is a complex $L^\cdot$ 
  of finitely generated projective $R$-modules such that for any $M\in 
  \Mod(R)$ and for any $i\ge 0$  we have 
\[ 
\Ext^i(\Fscr,\Gscr\otimes_R M)=H^i(L^\cdot \otimes_R M) 
\] 
\end{proposition} 
\begin{proof} 
\setcounter{step}{0} 
\begin{step} We first claim that there is an $N$ such that for $n\ge 
  N$ one has that $\Gamma(X,\Gscr(n))$ is a projective $R$-module, 
  $\Gamma(X,\Gscr(n)\otimes_R M)=\Gamma(X,\Gscr(n))\otimes_R M$ and 
  $R^i\Gamma(X,\Gscr(n)\otimes_R M)=0$ for $i>0$ and all $M$. We start 
  with the last part of this claim. We select $N$ is such a way that 
  $R^i\Gamma(X,\Gscr(n))=0$ for $i>0$ and $n\ge N$.  Using the fact that 
  $\Gamma(X,-)$ has finite cohomological dimension and degree shifting 
  in $M$ we deduce that indeed $R^i\Gamma(X,\Gscr(n)\otimes_R M)=0$ for 
  $i>0$ and all $M$. Thus $\Gamma(X,\Gscr(n)\otimes_R -)$ is an exact 
  functor. Applying this functor to a projective presentation of $M$ 
  yields $\Gamma(X,\Gscr(n)\otimes_R M)=\Gamma(X,\Gscr(n))\otimes_R M$. 
  Since $\Gamma(X,\Gscr(n)\otimes_R -)$ is left exact and 
  $\Gamma(X,\Gscr(n))\otimes -$ is right exact this implies that 
  $\Gamma(X,\Gscr(n))$ is flat. Finally since $A$ satisfies $\chi$ and 
  $R$ is noetherian $\Gamma(X,\Gscr(n))$ is finitely presented and hence 
  projective. 
\end{step} 
\begin{step} Now let $N$ be as in the previous step and take a resolution 
$\Pscr^\cdot \r \Fscr\r 0$ 
where the $\Pscr_i$ are finite direct sums of objects $\Oscr(-n)$ with 
$n\ge N$. 

Then $\Ext^i(\Fscr,\Gscr\otimes_R M)$ is the homology of 
\[ 
\Hom(\Pscr^\cdot,\Gscr\otimes_R M))=\Hom(\Pscr^\cdot,\Gscr)\otimes_R M 
\] 
where the equality follows from Step 1. We put 
$L^\cdot=\Hom(\Pscr^\cdot,\Gscr)$ which is term wise projective, also 
by Step 1. This finishes the proof. \qed 
\end{step} 
\def\qed{}\end{proof} 
For a point $x\in \Spec R$ we denote the base change functor 
$-\otimes_R k(x)$ by $(-)_x$. We also put $X_x=\Proj A_x$. 

\begin{corollary} 
\label{ref-A.3-64} 
If  $\Gscr$ is  as in the previous  proposition then the function 
\[ 
\Spec R\r \NN: x\mapsto \dim_{k(x)} R\Gamma^i(X_x,\Gscr_x) 
\] 
is upper semi-continuous. 
\end{corollary} 

\begin{proof}  By \cite[Lemma C6.6]{AZ1} 
\[ 
R\Gamma^i(X_x,\Gscr_x)=R\Gamma^i(X,\Gscr\otimes_R k(x)) 
\] 
This implies 
\begin{equation} 
\label{ref-A.1-65} 
R\Gamma^i(X_x,\Gscr_x)=H^i(L^\cdot\otimes_R k(x)) 
\end{equation} 
The fact that the dimension of the right hand side of 
\eqref{ref-A.1-65} is upper semi-continuous is an elementary fact 
from linear algebra. 
\end{proof} 
\begin{corollary} 
\label{ref-A.4-66} 
Assume that  $\Gscr$ is  as in the previous  proposition and assume that 
$R$ is a domain. Assume furthermore that the function 
\[ 
\Spec R\r \NN: x\mapsto \dim_{k(x)} R\Gamma^i(X_x,\Gscr_x) 
\] 
is constant. Then  $R\Gamma^i(X,\Gscr)$ is projective over $R$ and 
in addition for any $M\in \Mod(R)$ the natural map 
\[ 
R\Gamma^i(X,\Gscr)\otimes_R M\r  R\Gamma^i(X,\Gscr\otimes_R M) 
\] 
is an isomorphism for all $x\in  \Spec R$. 
\end{corollary} 
\begin{proof} This is proved as \cite[Corollary 12.9]{H}. 
\end{proof} 

\section{Hilbert series up to invariant $6$} \label{ref-B-67} 
\label{ref-B-68} 
\def\mystrut{\vrule width 0em height 1em} 
Let $I$ be a normalized rank one torsion free graded $A$-module of 
projective dimension one over a quantum 
polynomial ring with invariant $n$.  According to Theorem A the Hilbert 
series of $I$ has 
the form $h_{I}(t) = 1/(1-t)^{3} - s_I(t)/(1-t)$ where $s_I(t)$ is a 
Castelnuovo polynomial of weight $n$.  For the cases $n \leq 6$ we 
list the possible Hilbert series for $I$, the corresponding 
Castelnuovo polynomial, the dimension of the stratum (given by $\dim_{k} 
\Ext^{1}_A(I,I)$) and the possible minimal resolutions of 
$I$. 

{\small 
\[ 
\begin{array}{|c|l|} 
\hline 
n = 0 
& 
h_{I}(t)  = 1 + 3t + 6t^{2} + 10t^{3} + 15t^{4} + 21t^{5} + \ldots 
\mystrut\\ 
&s_{I}(t)  = 0 \\ 
& \dim_{k}  \Ext^{1}(I,I) = 0 \\ 
&0  \r A \r I \r 0\\ 
\hline 
n = 1& 
h_{I}(t)  = 2t + 5t^{2} + 9t^{3} + 14t^{4} + 20t^{5} + 27t^{6} + \ldots 
\mystrut\\ 
\parbox[c]{0.5cm}{\includegraphics[height=0.5cm,width=0.5cm]{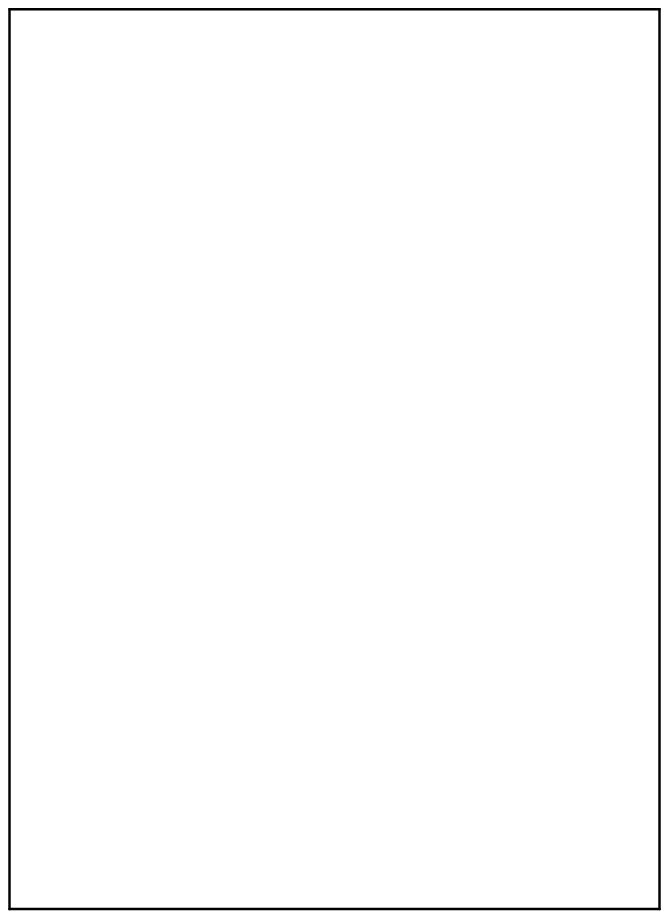}}& 
s_{I}(t)  = 1\\ 
&\dim_{k} \Ext^{1}(I,I) = 2 \\ 
&0  \r A(-2) \r A(-1)^{2} \r I \r 0\\ 
\hline 
n = 2& 
h_{I}(t)  = t + 4t^{2} + 8t^{3} + 13t^{4} + 19t^{5} + 26t^{6} + \ldots 
\mystrut\\ 
\parbox[c]{1cm}{\includegraphics[height=0.5cm,width=1cm]{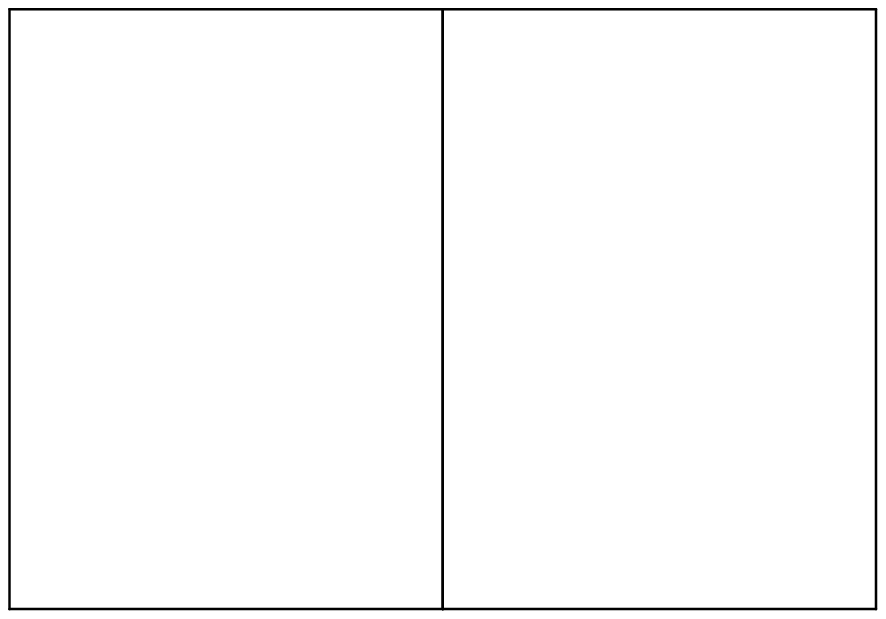}} 
&s_{I}(t)  = 1 + t 
 \\ 
&\dim_{k} \Ext^{1}(I,I) = 4 \\ 
&0  \r A(-3) \r A(-1) \oplus A(-2) \r I \r 0\\ 
\hline 
n = 3 
&h_{I}(t)  = 3t^{2} + 7t^{3} + 12t^{4} + 18t^{5} + 25t^{6} + \ldots \mystrut\\ 
\smash{\raisebox{-0.7cm}{\includegraphics[height=1cm,width=1cm]{fig3a.eps}}} 
&s_{I}(t)  = 1 + 2t 
  \\ 
&\dim_{k}  \Ext^{1}(I,I) = 6 \\ 
&0  \r A(-3)^{2} \r A(-2)^{3} \r I \r 0 \\ 
\cline{2-2} 
&h_{I}(t)  = t + 3t^{2} + 7t^{3} + 12t^{4} + 18t^{5} + 25t^{6} + \ldots 
\mystrut\\ 
\parbox[c]{1.5cm}{\includegraphics[height=0.5cm,width=1.5cm]{fig3b.eps}} 
&s_{I}(t)  = 1 + t + t^{2} 
\quad \\ 
&\dim_{k} \Ext^{1}(I,I) = 5 \\ 
&0  \r A(-4) \r A(-1) \oplus A(-3) \r I \r 0\\ 
\hline 
n = 4& 
h_{I}(t)  = 2t^{2} + 6t^{3} + 11t^{4} + 17t^{5} + 24t^{6} + \ldots \mystrut\\ 
\smash{\raisebox{-0.7cm}{\includegraphics[height=1cm,width=1.5cm]{fig4a.eps} 
}} &s_{I}(t)  = 1 + 2t + t^{2} 
\\ 
&\dim_{k} \Ext^{1}(I,I) = 8 \\ 
&0  \r A(-4) \r A(-2)^{2} \r I \r 0\\ 
&0  \r A(-3)\oplus A(-4) \r A(-2)^{2} \oplus A(-3)\r I \r 0\\ 
\cline{2-2} 
&h_{I}(t)  = t + 3t^{2} + 6t^{3} + 11t^{4} + 17t^{5} + 24t^{6} + \ldots 
\mystrut\\ 
\parbox[c]{2cm}{\includegraphics[height=0.5cm,width=2cm]{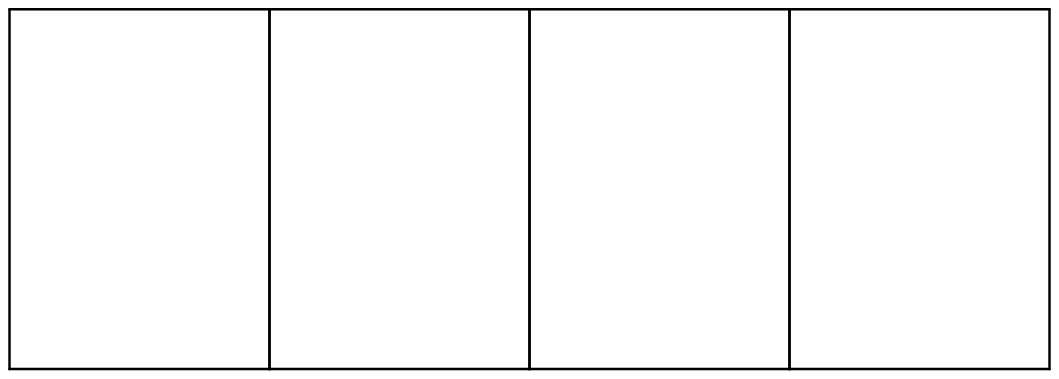}}& 
s_{I}(t)  = 1 + t + t^{2} + t^{3}\\ 
&\dim_{k} \Ext^{1}(I,I) = 6 \\ 
&0  \r A(-5) \r A(-1) \oplus A(-4) \r I \r 0\\ 
\hline 
n = 5 
&h_{I}(t)  = t^{2} + 5t^{3} + 10t^{4} + 16t^{5} + 23t^{6} + \ldots \mystrut\\ 
\smash{\raisebox{-0.7cm}{\includegraphics[height=1cm,width=1.5cm]{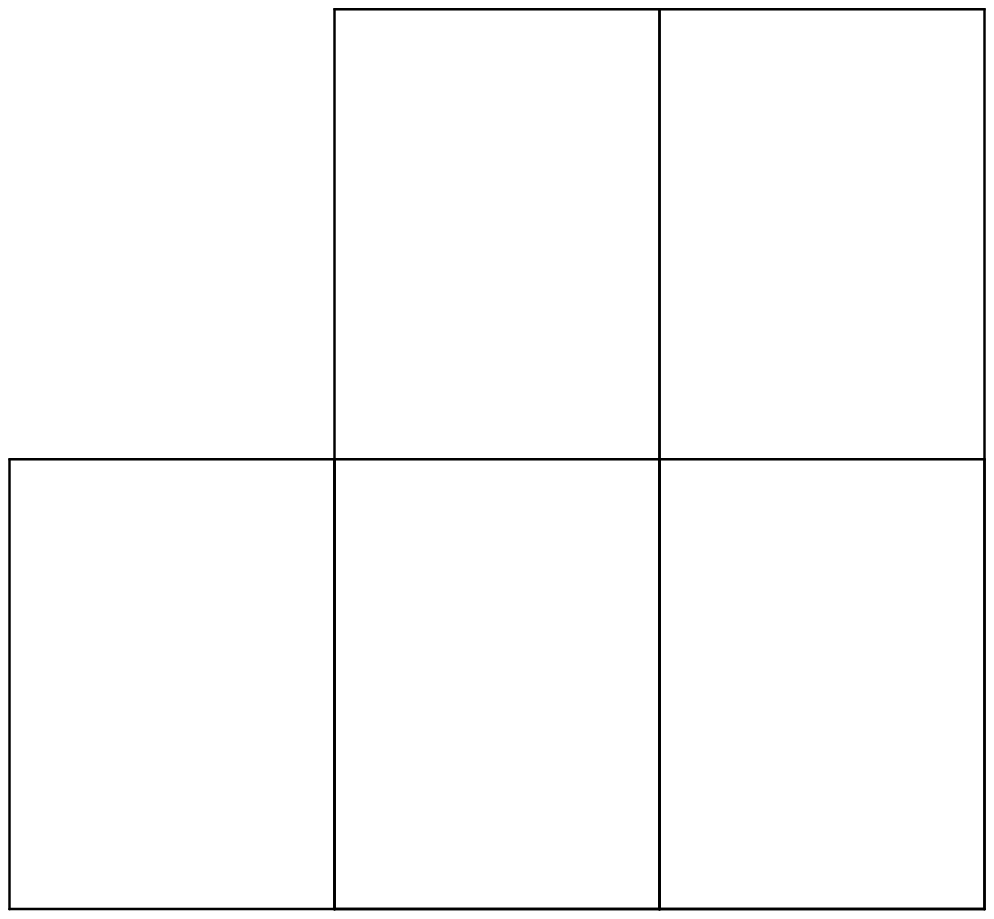} 
}}&s_{I}(t)  = 1 + 2t + 2t^{2} 
\\ 
&\dim_{k} \Ext_A^{1}(I,I) = 10 \\ 
&0  \r A(-4)^{2} \r A(-2) \oplus A(-3)^{2} \r I \r 0\\ 
\cline{2-2} 
&h_{I}(t)  = 2t^{2} + 5t^{3} + 10t^{4} + 16t^{5} + 23t^{6} + \ldots \mystrut\\ 
\smash{\raisebox{-0.7cm}{\includegraphics[height=1cm,width=2cm]{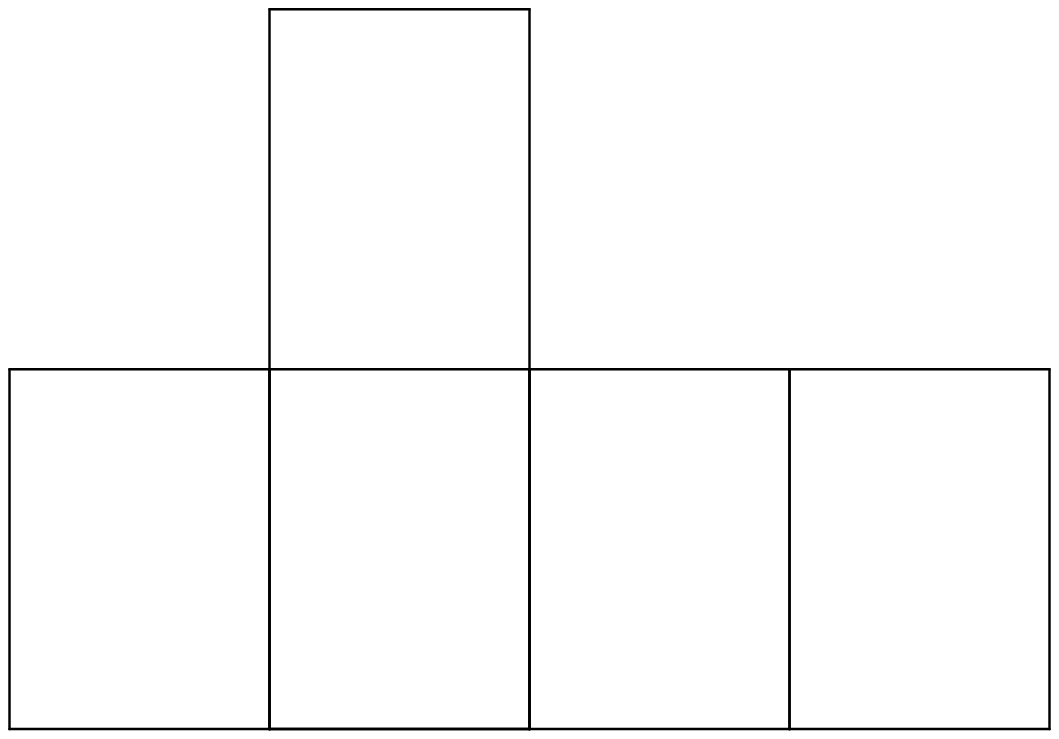}}}& 
s_{I}(t) = 1 + 2t + t^{2} + t^{3} 
\\ 
&\dim_{k} \Ext_A^{1}(I,I) = 8 \\ 
&0  \r A(-3) \oplus A(-5) \r A(-2)^{2} \oplus A(-4) \r I \r 0\\ 
\cline{2-2} 
&h_{I}(t)  = t + 3t^{2} + 6t^{3} + 10t^{4} + 16t^{5} + 23t^{6} + 
\ldots 
\mystrut\\ 
\smash{\parbox[c]{2.5cm}{\includegraphics[height=0.5cm,width=2.5cm]{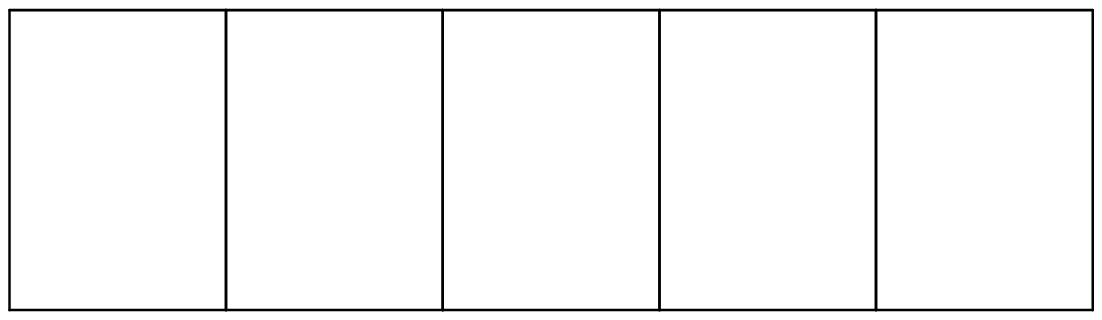}}} &s_{I}(t)  = 1 + t + t^{2} + t^{3} + t^{4} 
 \\ 
&\dim_{k} \Ext_A^{1}(I,I) = 7 \\ 
&0  \r A(-6) \r A(-1) \oplus A(-5) \r I \r 0\\ 
\hline 
\end{array} 
\] 
\newpage 
\[ 
\begin{array}{|c|l|} 
\hline 
n = 6 
&h_{I}(t)   = 4t^{3} + 9t^{4} + 15t^{5} + 22t^{6} + 30t^{7} + \ldots \mystrut\\ 
\smash{\raisebox{-1.1cm}{ 
\includegraphics[height=1.5cm,width=1.5cm]{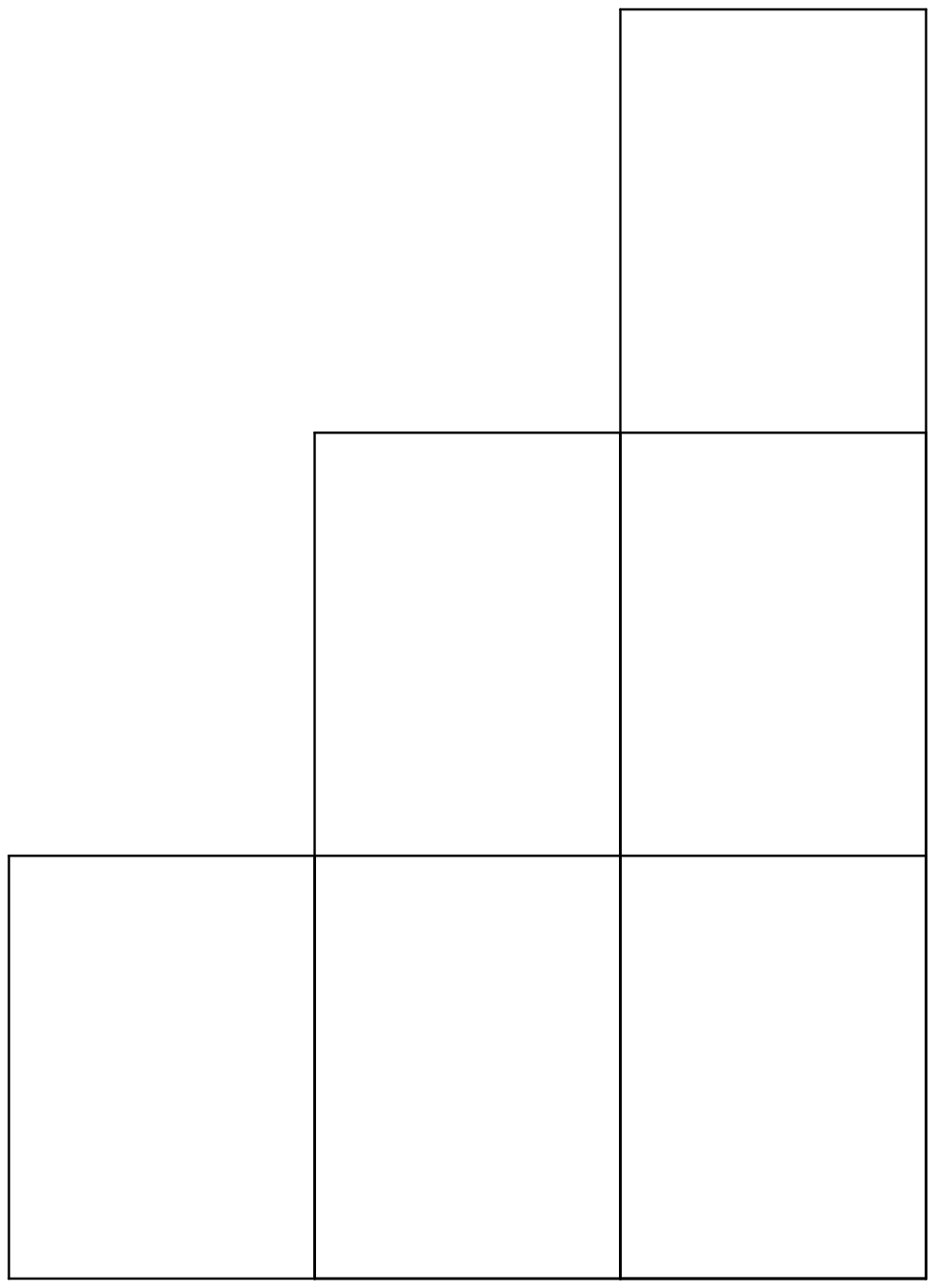}}} &s_{I}(t)  = 1 + 2t 
+ 3t^{2} 
 \\ 
&\dim_{k} \Ext_A^{1}(I,I) = 12 \\ 
&0  \r A(-4)^{3} \r A(-3)^{4} \r I \r 0\\ 
\cline{2-2} 
&h_{I}(t)  = t^{2} + 4t^{3} + 9t^{4} + 15t^{5} + 22t^{6} + \ldots\mystrut \\ 
\smash{\raisebox{-0.7cm}{\includegraphics[height=1cm,width=2cm]{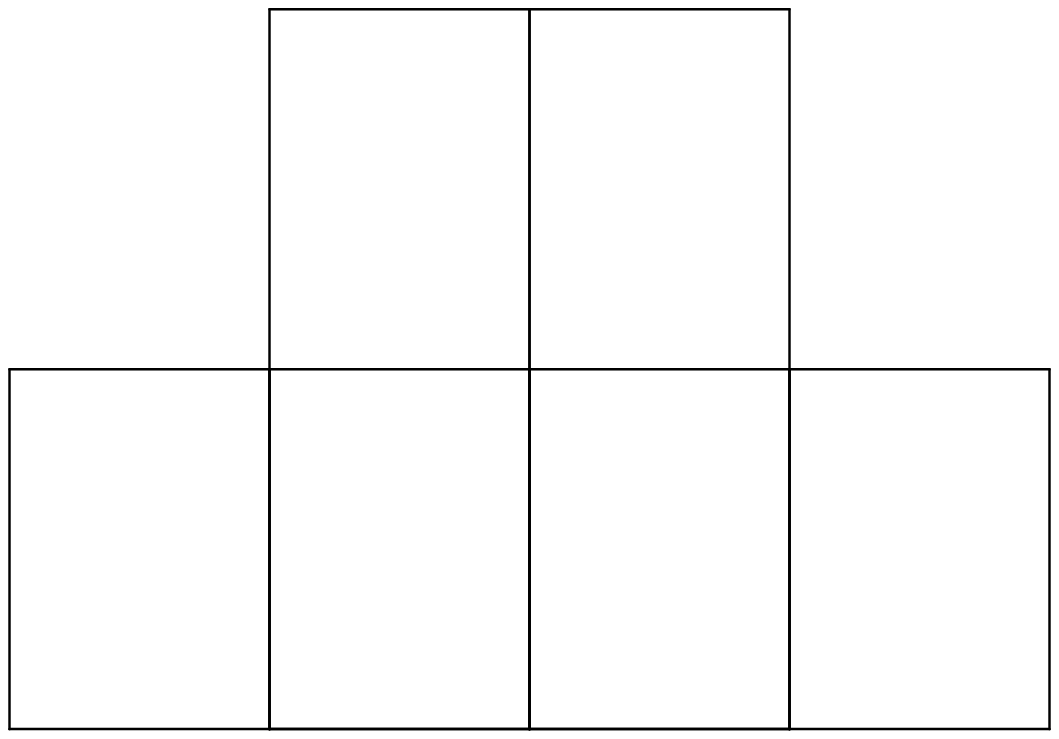}}}& 
s_{I}(t)  = 1 + 2t + 2t^{2} + t^{3}\\ 
&\dim_{k} \Ext_A^{1}(I,I) = 11 \\ 
&0  \r A(-5) \r A(-2) \oplus A(-3) \r I \r 0\\ 
&0  \r A(-4)\oplus A(-5) \r A(-2) \oplus A(-3)\oplus A(-4) \r I \r 0\\ 
\cline{2-2} 
&h_{I}(t)  = 2t^{2} + 5t^{3} + 9t^{4} + 15t^{5} + 22t^{6} + \ldots \mystrut\\ 
\smash{\raisebox{-0.7cm}{\includegraphics[height=1cm,width=2.5cm]{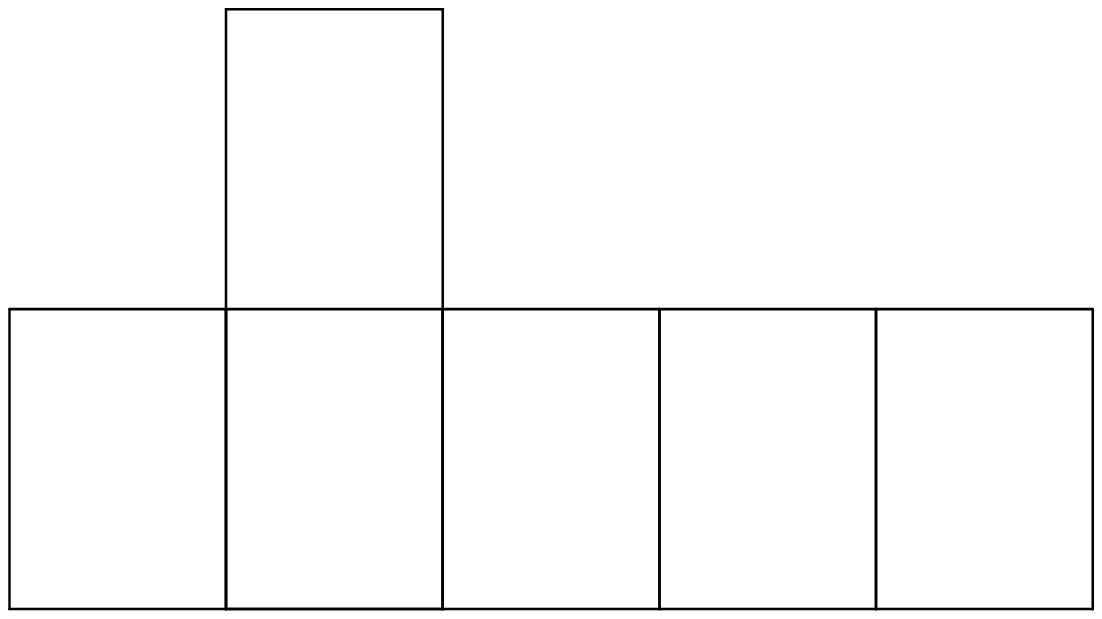} 
}}&s_{I}(t)  = 1 + 2t + t^{2} + t^{3} + t^{4}\\ 
&\dim_{k} \Ext_A^{1}(I,I) = 9 \\ 
&0  \r A(-3) \oplus A(-6) \r A(-2)^{2} \oplus A(-5) \r I \r 0\\ 
\cline{2-2} 
&h_{I}(t)  = t + 3t^{2} + 6t^{3} + 10t^{4} + 15t^{5} + 22t^{6} + \ldots 
\mystrut\\ 
\smash{\parbox[c]{3cm}{\includegraphics[height=0.5cm,width=3cm]{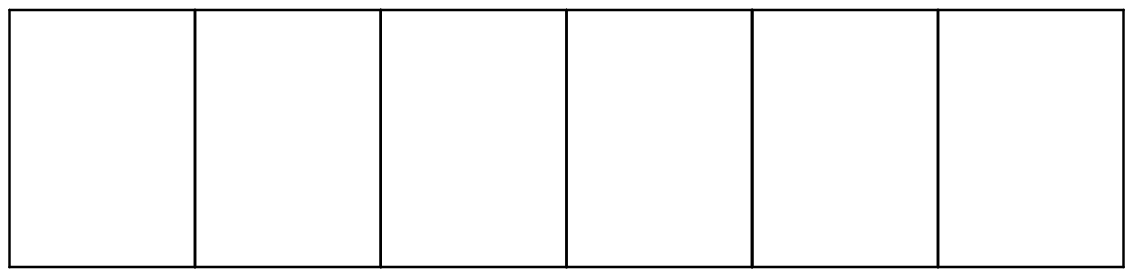}}} 
&s_{I}(t)  = 1 + t + t^{2} + t^{3} + t^{4} + t^{5} 
 \\ 
&\dim_{k} \Ext_A^{1}(I,I) = 8 \\ 
&0  \r A(-7) \r A(-1) \oplus A(-6) \r I \r 0\\ 
\hline 
\end{array} 
\] 
} 
\ifx\undefined\bysame 
\newcommand{\bysame}{\leavevmode\hbox to3em{\hrulefill}\,} 
\fi

\end{document}